\documentclass[a4paper,11pt,reqno]{amsart}
\usepackage{epsfig,url,paralist}
\usepackage{amssymb}
\usepackage[normalem]{ulem}
\usepackage{amsmath}
\usepackage{amsthm, tikz}
\usepackage{hyperref,verbatim}
\usepackage{graphicx}
\usepackage{setspace}
\usepackage{enumitem,lineno}
\usepackage{calrsfs}
\usepackage{soul} 
\usepackage{array, multirow, bigdelim, makecell, booktabs} 
\usepackage{cleveref}
\numberwithin{equation}{section}

\newcommand{\hvm}[1]{\textcolor{red}{#1}}

\newtheorem{theorem}{Theorem}
\newtheorem*{theorem*}{Main Theorem}
\newtheorem*{conjecture*}{Wiegold conjecture}

\newtheorem{lem}[theorem]{Lemma}
\newtheorem{obs}[theorem]{Observation}

\newtheorem{prop}[theorem]{Proposition}

\numberwithin{theorem}{section}
\newcommand{\PSL}{\mathsf{PSL}}

\newcommand{\Aut}{\mathsf{Aut}}
\newcommand{\Syl}{\mathsf{Syl}}

\newcommand{\PGL}{\mathsf{PGL}}

\newcommand{\SL}{\mathsf{SL}}

\newcommand{\PG}{\mathsf{PG}}

\newcommand{\U}{\mathsf{U}}
\newcommand{\Alt}{\mathsf{Alt}}

\def\<{\langle}
\def\>{\rangle}

\usepackage{fancyhdr}
\begin{document}
\title{On Wiegold's conjecture for the small Ree groups}
\author[Busch, Pengitore, Schillewaert, Van Maldeghem]{Sira Busch, Mark Pengitore, Jeroen Schillewaert and Hendrik Van Maldeghem}
\thanks{The research of all authors is supported by the New Zealand Marsden Fund. The second author is supported by the National Science Center Grant Maestro-13 UMO-2021/42/A/STI/00306}
\thanks{The fourth author is partially supported by Ghent University Special Research Fund, grant BOF.24Y.2021.0029.01}
\date{}
\maketitle

\begin{abstract}
The Wiegold conjecture holds for the small Ree groups for $k$-tuples where $k \geq 5$.
\end{abstract}
\section{Introduction}

We prove the following for the small Ree groups:
\begin{theorem*} \label{main_thm1}
The small Ree groups $^2 \mathsf{G}_2(3^{2e+1})$ satisfy the Wiegold conjecture for $k$-tuples, for all $e\geq 1$ and all $k \geq 5$. 
\end{theorem*}

The original motivation of the Wiegold conjecture comes from actions of the automorphism group of the free group on group presentations.

\subsection{$T_k$-systems and the Wiegold conjecture}

Let $G$ be a group. We call $N$ a $G$-defining subgroup of a free group $F_k$ if $N\leq F_k$ and $F_k/N\cong G$. If we denote by
$\Sigma_k(G)$ the set of $G$-defining subgroups of $F_k$ we can consider an action of $\Aut(F_k)$ on $\Sigma_k(G)$. The orbits of this action are $T_k$ systems introduced by Bernhard and Hanna Neumann \cite{NN}. The group ${\Aut}(F_k)$ is generated by Nielsen moves, transpositions and inversions (which we will all call Nielsen moves henceforth). These moves thus define a graph structure on $\Sigma_k(G)$. Thus, $T_k$-systems correspond to connected components of $\Sigma_k(G)$.

The unpublished conjecture below, likely inspired by Gilman \cite{Gilman}, is due to Wiegold. Graham and Diaconis \cite{DG} made a similar conjecture for $S_n$.
\begin{conjecture*}
For every $k\geq 3$ and every nonabelian simple group $G$, there exists only one $T_k$-system, i.e. the graph $\Sigma_k(G)$ is connected.
\end{conjecture*}

Thus far, the following was known on the above conjecture:
\begin{theorem}
The Wiegold conjecture holds for $k$-tuples in these cases.
\begin{itemize}
\item \cite{Gilman} $G = \PSL(2,p)$, where $p\geq 5$ and $k\geq 3$
\item \cite{E93} $G = \PSL(2,2^m)$ or $\mathsf{Sz}(2^{2m-1})$, where $m\geq 2$ and $k\geq 3$
\item \cite{G08} $G = \PSL(2,q)$ where $q$ is odd and $k\geq 4$.
\item \cite{David} $G = \Alt(6),\Alt(7)$ where $k=3$.
\item \cite{Pak} $G = \Alt(8),\Alt(9),\Alt({10})$ where $k=3$.
\end{itemize}
\end{theorem}

Avni and Garion \cite{Avni-Garion} used deep results of Larsen and Pink \cite{Larsen-Pink} on finite subgroups of algebraic groups to prove that for finite simple groups of Lie type in characteristic different from $2$ and $3$ there exists a function $c$, only depending on the Lie rank $r$, such that the Wiegold conjecture holds for any $k \geq c(r)$. This is independent of the field but grows exponentially in $r$. Avni informed us that even if their methods would be made to work for Ree groups the constant would be huge.

\subsection{Connections of the Wiegold conjecture}

Lubotzky \cite{Lubotzky} wrote a very interesting survey on the dynamics of $\Aut(F_n)$ actions which is closely related to the Wiegold conjecture. Some connections include the following:
\begin{itemize}
\item Nielsen equivalence classes of generators of a group $G$ parameterise free actions of $G$ on handlebodies \cite[Theorem 3.4]{Lubotzky}.
\item The Wiegold conjecture is equivalent to: for $n\geq 3$ and G a finite simple group the action of ${\Aut}(F_n)$ on $\mathsf{Epi}(F_n,G)$ is transitive. In \cite{Gelander} Gelander considered an analogous question for compact Lie groups, non-compact simple analytic groups and simple algebraic groups.
\item Product replacement algorithm which we discuss in detail below.
\end{itemize}

\subsection{Product replacement algorithm graph}

The product replacement algorithm was designed by Leedham-Green and Soicher as mentioned in the first description of it by Celler, Leedham-Green, Murray, Niemeyer and O'Brien \cite{EOB}. It very efficiently generates random elements in a finite group. Lubotzky and Pak \cite{LubotzkyPak} explained why the product replacement algorithm works so well using Kazhdan's property (T). We closely follow Igor Pak's survey paper \cite{Pak}.

Let $n \geq 2$ and let $G$ be a finite simple group. Let
\[
V_n(G) = \{(g_1, \ldots, g_n) \in G^n : \left<g_1, \ldots, g_n \right> = G\}.
\]
{The set $V_n(G)$ is the set of all generating $n$-tuples of $G$.} Such an $n$-tuple is called \emph{minimal} if no proper subtuple generates $G$. Otherwise it is called \emph{redundant}. 

The \emph{extended PRA graph} $\widetilde{X}_n(G)$ has vertices given by $V_n(G)$ and edges corresponding to the \emph{Nielsen moves} $R_{i,j}^{\pm}, L_{i,j}^{\pm}, P_{i,j}, I_i,$ for $1 \leq i \neq j \leq n,$ where
\[
R_{i,j}^{\pm} \colon (g_1, \ldots, g_i, \ldots, g_n) \to (g_1, \ldots, g_i g_j^{\pm 1}, \ldots, g_n),
\]
\[
    L_{i,j}^{\pm} \colon (g_1, \ldots, g_i, \ldots, g_n) \to (g_1, \ldots,  g_j^{\pm 1}g_i, \ldots, g_n),
\]
\[
P_{i,j} \colon (g_1, \ldots, g_i, \ldots, g_j, \ldots, g_n) \to  (g_1, \ldots, g_j, \ldots, g_i, \ldots, g_n),
\]
and
\[
I_i \colon (g_1, \ldots, g_i, \ldots, g_n) \to (g_1, \ldots, g_i^{-1}, \ldots, g_n).
\]

The Wiegold conjecture can be reformulated as:
\begin{conjecture*}
Let $G$ be a nonabelian finite simple group. Then $\widetilde{X}_n(G)$ is connected for every $n \geq 3.$
\end{conjecture*}

By a celebrated recent result of Burness, Guralnick and Harper \cite{BGH}, finite simple groups have spread at least two and hence, 
the graph induced on the redundant $n$-tuples of $\widetilde{X}_n(G)$ is connected, see e.g. \cite[Proposition 2.5.12]{Pak}. Hence, in order to show that $\widetilde{X}_n(G)$ is connected, it suffices to show that every vertex of $\widetilde{X}_n(G)$ can be connected via a path to a redundant vector. In the rest of the paper, we will use the phrase \emph{connected to a vector} to mean that two given vectors are in the same connected component of $\widetilde{X}_n(G)$ (and we will only consider $n\in\{3,4,5\}$ and $G$ a small Ree group.  Note that this means that we can change a given vector using only Nielsen moves into another given vector. 

\subsection{Wiegold and its variations for more general groups}

Let $d(G)$ denote the minimum number of generators of $G$. Pak
\cite{Pak} asked whether there are finite groups such that the extended PRA graph  $\widetilde{X}_k(G)$ is disconnected for $k\geq d(G)+1$. As $d(G)=2$ for finite simple groups, conjecturing this does not happen is a generalisation of the Wiegold conjecture.
This generalisation was proved for abelian groups by the Neumann and Neumann \cite{NN} and for solvable groups by Dunwoody \cite{DW}. Garion \cite{G08} also proved it for $\PGL(2,q)$ where $q$ is an odd prime power. 
The next natural class of (almost) simple groups to look at are the projective special (general) unitary groups, which would conclude the rank one case of the Wiegold conjecture. Although a number of arguments here translate directly to this case extra complications arise from the richer subgroup structure and the fact that the characteristic is not uniform. 

More broadly, Nielsen equivalence classes of generating tuples have also been studied extensively for infinite groups. For example Kapovich and Weidmann \cite{KW} provided a counterexample to a Wiegold type conjecture for word hyperbolic groups. 
If one wishes to disprove (variations of) the Wiegold conjecture for a given group one needs to show that certain generating $n$-tuples are not Nielsen equivalent. For concrete small groups this can, in principle, be done by exhaustive computation.
For larger groups or families of groups a potential strategy consists in using Fox calculus which can provide obstructions to Nielsen equivalence \cite{Lus91,LM93}.

As pointed out by Kapovich and Weidmann \cite{KW}, Nielsen equivalence is not decidable for finitely presented torsion-free small cancellation groups since they do not have a decidable subgroup membership problem (which is a special case) as shown by Rips \cite{Rips}.

\subsection{Proof outline}
We basically follow the same strategy as Evans \cite{E93} introduced to handle the Suzuki groups and as was also used later by Garion \cite{G08} for $\PSL_2(q)$. So, we start off with a generating $k$-tuple of $G:={^2\mathsf{G}}_2(q)$, $k\in\{3,4,5\}$ (we do as much as we can for minimal $k\in\{3,4,5\}$, the proofs for larger $k$ being totally similar to the ones for smaller values of $k$). The eventual goal is to show that it is connected to a redundant vector. To accomplish that, there are three major steps. They all use the natural $2$-transitive permutation representation of $G$ on a set $U_R(q)$ with $q^3+1$ points. Unlike the situation for the Suzuki groups and $\PSL_2(q)$, this set is in fact a rank 2 geometry (and not rank 1) owing to the existence of simple subgroups of Lie  type over $\mathbb{F}_q$, namely $\PSL_2(q)$. The latter subgroups stabilise sets of $q+1$ points on which they act in the natural way (and so we can think of each such set of $q+1$ points as a projective line over $\mathbb{F}_q$; we call it a block). The set of blocks turns $U_R(q)$ into a unital, that is, a $2-(q^3+1,q+1,1)$ design \cite{Onan}. Now, the eventual goal is to show that we can connect to a $k$-tuple containing an element $x$ with the property that all $(k-1)$-subtuples containing $x$ generate subgroups contained in a subfield subgroup and in no maximal subgroup of another type, and the same holds for the subgroups generated by $x$ and the conjugates of $x$ by any $k-2$ other elements of the vector. The first step to accomplish this consists in showing that we can connect to a vector containing an element $x$ that is neither unipotent nor involutive (unipotent elements and involutions are precisely the elements of the Ree group that either fix exactly 1 or $q+1$ points of the corresponding Ree unital). In this first step, understanding the structure of the point stabiliser is essential and for that, we use the explicit description using $7\times7$ matrices over $\mathbb{F}_q$ provided in \cite{TTM07}. We survey this description in the appendix, where we also prove some more details about it. In the second step, accomplished in Sections~\ref{sec:structural} and~\ref{sec:redundant-connect}, we first show in \ref{sec:structural} that we can connect to a $k$-tuple such that the above subgroups are not contained in either a point-stabiliser or a block stabiliser. Then in the first part of \Cref{sec:redundant-connect} we further eliminate subgroups contained in  maximal subgroups other than a point-stabiliser, a block stabiliser or a subfield subgroup. In Step 3, accomplished in the remainder of \ref{sec:redundant-connect}, we use an inductive procedure (inspired by Evans \cite{E93})  to connect to a redundant vector.

Each of these steps needs some additional ideas, compared to \cite{E93,G08}, to make the strategy work. The Ree groups are generally known to have a rather peculiar structure. For instance,  Guralnick, Kantor, Kassabov and Lubotzky \cite{GKKL} provide short and small presentations for all finite simple groups---except for the Ree groups, the latter has been done recently by Hulpke, Kassabov, Seress and Wilson \cite{HKSW}.

On the positive side, the list of maximal subgroups of $^2\mathsf{G}_2(q)$ is rather small and can be compared to the Suzuki groups. 

{\bf Convention.} Throughout the paper lemmas and proofs might be formulated for $k$-tuples, where in fact the corresponding statements and proofs are valid for $\ell$-tuples, $\ell\geq k$. We opted to do this for clarity of exposition; the reader can easily check there is no harm in taking longer tuples as one can just leave the appropriate $\ell-k$ entries unchanged when repeating the appropriate arguments for a different $k$-sub-tuple.

{\bf Acknowledgement.} We thank Alex Lubotzky who generously shared his insights on the Wiegold conjecture and pointed out various references.

\section{The small Ree groups}

\subsection{General notation}

Given a finite group $G$ and a prime $p$, we denote $\Syl_p(G)$ as the collection of Sylow $p$-subgroups. For a natural number $m$, we denote the cyclic subgroup of order $m$ by $C_m$. For a subgroup $H \leq G$, we denote the centraliser of $H$ in $G$ as $C_G(H)$ and the normaliser of $H$ in $G$ as $N_G(H)$. We denote the commutator subgroup of $G$ by $[G,G]$, by $Z(G)$ the center of $G$, and by $\Phi(G)$ the Frattini subgoup of $G$. Given a group, we denote $|G|$ as its order and $|x|$ as the order of element $x \in G$. We denote the greatest common divisor of two integers $n$ and $m$ as $\gcd(m,n)$.

A \emph{Frobenius group} $G$ is a transitive permutation group on a finite set $X$, such that no nontrivial element fixes more than one point and some nontrivial element fixes a point. A subgroup $H \leq G$ fixing a point of $X$ is called a \emph{Frobenius complement}. The identity element together with all elements not in any conjugate of $H$ form a normal subgroup called the \emph{Frobenius kernel $K$}. In particular, $G = K \rtimes H$. There is a unique Frobenius group of order $21$ and we denote it by $\mathsf{Frob}(21)$. Given a finite group $G$, a subgroup $H \leq G$ is called a \emph{Hall subgroup} if $\gcd(|H|, [G:H]) = 1$.

Finally we consider actions on the left; hence, our convention is that $x^y=yxy^{-1}$ and $^yx=y^{-1}xy$.

\subsection{Generalities}\label{sec:generalities}
There are three types of Suzuki--Ree groups: type ${^2\mathsf{B}}_2$, usually called \emph{Suzuki groups}, types $^2\mathsf{G}_2$ and $^2\mathsf{F}_4$, usually called \emph{small} and \emph{large Ree groups}, respectively. Henceforth, \emph{Ree group} will mean one of type $^2\mathsf{G}_2$.
The Ree group $G={^2\mathsf{G}}_2(q)$, sometimes also denoted as $\mathsf{R}(q)$, has order $q^3(q+1)(q-1)$ where $q = 3^{2e+1}$. We denote by $\theta$ the field automorphism $x\mapsto x^{3^{e+1}}$; note $x^{\theta^2}=x^3$. Hence, $\theta$ is the square root of the Frobenius automorphism. The subgroup structure of the Ree groups was described by Levchuk and Nuzhin \cite{LN85}, see also Kleidman's paper \cite{K88} and  \cite[Section~9.2.4]{TTM07}. The Ree groups are one of the families of exceptional finite simple groups of Lie type. From now on we fix $q=3^{2e+1}$, $e\geq 1$, and denote the Ree group $^2\mathsf{G}_2(q)$ by $G$. Note that $e=0$ corresponds to $\PSL_2(8)\rtimes C_3\cong \mathsf{P\Gamma L}_2(8)$.

The Ree group $G$ has a natural doubly transitive action on its set $\Omega:=\Syl_3(G)$ of Sylow $3$-subgroups. The stabiliser $G_P$ of $P\in \Omega$ is given by $N_G(P)$. The stabiliser $G_{P, Q}$ of two members of $\Syl_3(G)$ is a cyclic subgroup $C_{q-1}$.  

A Sylow $2$-subgroup $P$ of $G$ is elementary abelian of order $8$. Moreover $C_G(P) = P$ and $N_G(P) \cong 2^3\rtimes\mathsf{Frob}(21)$, which has order $8 \cdot 7 \cdot 3=168$. As a result, the $2$-subgroups of equal order are conjugate. In particular, all involutions are conjugate to one another. The centraliser of an involution $\eta$ is isomorphic to $\PSL_2(q)\times\<\eta\>$.

The group $G$ has cyclic Hall subgroups $M_i$ ($i = \pm 1)$ of orders $q+1+i \cdot \sqrt{3q}$. Each subgroup $M_i$ coincides with $C_G(x)$ for all nontrivial $x \in M_i$. The subgroups $N_G(M_i)$ are Frobenius groups with kernel $M_i$ and cyclic complement of order $6$. For each subgroup $V$ of order $4$, there exists a cyclic Hall subgroup $M_0$ of order $\frac{q+1}{4}$ and an element $t$ of order $6$ such that $N_G(V) = N_G(M_0) = V \rtimes (M_0 \rtimes \left<t\right>) \cong (V \times (M_0 \rtimes C_2)) \rtimes C_3$, where $M_0\rtimes C_2$ is dihedral of order $\frac{q+1}2$. All subgroups of order $\frac{q-1}{2}$, $\frac{q+1}{4}$, $q+\sqrt{3q} + 1$, or $q-\sqrt{3q}+1$ are conjugate in $G$. The same holds for subgroups of order $q-1$ and $q+1$.

A Sylow $3$-subgroup $P$ of $G$ has order $q^3$. The following properties can also be checked directly on the model presented in \Cref{appendixA} (but are well known). The center $Z(P)$ of $P$ is an abelian subgroup of order $q$, $P$ has nilpotent step length $3$, and $[P,P] = \Phi(P)$ is an elementary abelian subgroup of order $q^2$ containing $Z(P)$. The elements of $P \backslash [P,P]$ have order $9$ where their cubes are forming $Z(P) \backslash \{1\}$. The involutions of $N_{G}(P)$ are conjugate in $N_{G}(P)$. We collect some of the above information in a lemma. 

\begin{lem}
	The maximal subgroups of $G$ are exhausted up to conjugacy, by the following
	\begin{itemize}
		\item $N_G(P)\cong P\rtimes C_{q-1}$, the normaliser of a Sylow $3$-subgroup $P$;
		\item $C_G(\eta)\cong \PSL_2(q)\times\<\eta\> $, the centraliser of an involution $\eta$;
		\item $^2\mathsf{G}_2(q_0), q = q_0^\delta$, $\delta$ being a prime;
		\item $N_G(M_i)$, the normaliser of a cyclic Hall subgroup $M_i$, $i=-1,0,1$, of order $q-\sqrt{3q}+1$, $\frac{q+1}{4}$ and 
		$q+\sqrt{3q}+1$, respectively.	
	\end{itemize}	
\end{lem}

\subsection{The Ree unital $U_R(q)$} As mentioned earlier, $G$ has a natural $2$-transitive action on its set $\Omega$ of Sylow $3$-subgroups. We can give $\Omega$ the structure of a \emph{unital}, that is, a $2-(q^3+1,q+1,1)$ design, by defining the blocks to be the sets of members of $\Omega$ fixed by a given involution. No block is (pointwise) fixed by two distinct involutions; hence, there are \[\frac{|G|}{2\cdot|\PSL_2(q)|}=q^2(q^2-q+1)\] blocks. Every pair of Sylow $3$-subgroups is contained in a unique block, which we sometimes call the \emph{join} of the two corresponding elements of $\Omega$. In the sequel we will refer to the elements of $\Omega$ as \emph{points}, and we will denote the points by Greek lower case letters. The unital will be denoted by $U_R(q)$, and consists of points and blocks. The unique block containing two given distinct points $\alpha$ and $\beta$ will be denoted as $B(\alpha,\beta)$. 

If $\alpha\in\Omega$, it is in fact a Sylow $3$-subgroup, but we see it merely as a point. When we want to emphasise the group structure of $\alpha$, we sometimes denote it also as $P_\alpha$. 

We list some properties of the permutation group $(G,U_R(q))$. 
For a point $\alpha$ we denote the set of blocks through $\alpha$ by $\mathcal{B}_\alpha$. We say that an element $x\in G$ acts freely on $U_R(q)$ if $\langle x \rangle$ does.

\begin{lem}\label{actionofP}
Let $P\in\Syl_3(G)$. Then $P$ fixes a unique point $\alpha\in\Omega$ and acts sharply transitively on $\Omega\setminus\{\alpha\}$. In particular, it acts transitively on $\mathcal{B}_\alpha$.  Moreover, each element of $P\setminus[P,P]$, which has order $9$, acts freely on the set of blocks, in particular on $\mathcal{B}_\alpha$. Also, each element of $[P,P]\setminus Z(P)$ stabilises exactly $q^2$ blocks, all belonging to $\mathcal{B}_\alpha$, and acts freely on the set of all other blocks. Finally, each member of $Z(P)$, which is the third power of an element of order $9$, acts freely on the set of blocks. 
\end{lem}

\begin{proof}
The first two assertions are easy and immediate, noting $P=P_\alpha$. Now let $x\in P$ have order $9$. Suppose $x^i$ stabilises a block $B$, $i\in\{1,2,\ldots,8\}$. Since $|x^i|\in\{3,9\}$ and $|B|=q+1$, $x^i$ fixes a point on $B$, and hence, $\alpha\in B$ (by the first assertion). Also, $x$ acts freely on the set of blocks through $\alpha$ if, and only if, $x^3$ does. Now we use the notation of \Cref{appendixA}. Up to conjugacy, we may assume $x^3=(0,0,1)_\infty$. Then one calculates $x^3\cdot B_{b,b''}=B_{b,b''+1}$, for all $b,b''\in\mathbb{F}_q$. It follows that both $x$ and $x^3$ act freely on the set of blocks through $\alpha$. Now let $x=(0,a',0)_\infty$ be, up to conjugacy, an arbitrary member of $[P,P]\setminus Z(P)$. Then $x$ maps the block $B_{b,b''}$ to the block $B_{b,b''-a'b}$. It follows that the $q$ blocks $B_{0,b''}$, $b''\in\mathbb{F}_q$, are stabilised, and none of the other blocks is.  A similar argument as for $|x|=9$ shows that every element of $P$ of order $3$ acts freely on the set of blocks not containing $\alpha$. 

The lemma is proved.
\end{proof}

We need one more property of the Sylow $3$-subgroups.

\begin{lem}\label{actionofP'}
Let $P\in\Syl_3(G)$, and let $\alpha\in\Omega$ be fixed by $P$. Then the set of all subsets of $\mathcal{B}_\alpha$ that are the sets of fixed lines of members of $P$, forms a system of imprimitivity for $P$, pointwise fixed under the action of $[P,P]$. Also, no element of $P\setminus[P,P]$ stabilises any member of that partition; hence, every element of $P\setminus[P,P]$ has orbits of size $3$. 
\end{lem}

\begin{proof}
We again use the notation of \Cref{appendixA} and put $\alpha=\infty$. Using the explicit form of $[P_\infty,P_\infty]$, where $P=P_\infty$, given there, it is easy to see that that the system of imprimitivity is \[\{B_{a,a''}\mid a''\in\mathbb{F}_q\}\mid a\in\mathbb{F}_q\}.\] The lemma is now obvious from the fact that each member of $P\setminus[P,P]$ has the form $(b,b',b'')_\infty$ with $b\neq 0$. 
\end{proof}

\color{black}
\begin{lem}\label{lem:no-blocks}
An element $x$ fixing at least three points not contained in a common block is the identity. An element $x$ fixing exactly two points $\alpha$ and $\beta$ stabilises a unique block, namely $B(\alpha,\beta)$.
\end{lem}
\begin{proof}
The first assertion follows straight from \Cref{obs1}. 
Now assume that an element $x\in G_{\alpha,\beta}$ stabilises a second block $B\neq B(\alpha,\beta)$. By the second assertion of \Cref{obs1}, $B\setminus B(\alpha,\beta)$ must contain a multiple of $|x|$ points. Since $|x|$ divides $q-1$, it can only divide also $|B\setminus B(\alpha,\beta)|$ if $|x|=2$ and $B\cap B(\alpha,\beta)=\varnothing$.  But then $x$ is an involution and fixes $B(\alpha,\beta)$ pointwise, contrary to our assumption that $x$ fixes exactly two points. 
\end{proof}

{\bf Orders and fixpoints} We will use the following fixpoint properties on a few occasions. Although they must be well-known, we sketch short proofs. 
\begin{lem}\label{lem:order-fix}Let $x\in G$ have order $k>1$. The one of the following holds.
\begin{itemize}
\item $k=2$ and the points fixed by $x$ are the $q+1$ points of a block.
\item $k\in\{3,6,9\}$ and $x$ fixes a unique point.
\item $k|\frac{q-1}2$ and $x$ fixes exactly two points.
\item $k|\frac{q+1}4$ and $x$ acts freely on $\Omega$, hence does not have any fixed points. 
\item $k|q^2-q+1$ and $x$ acts freely on $\Omega$, in particular it has no fixed points. 
\item $k\neq 2$ is even and $k|q-1$, then $x$ fixes exactly two points.
\item $k\neq 2$ is even and $k|\frac{q+1}2$, then $x$ has no fixed points.
\end{itemize}
\end{lem}

\begin{proof}
If $x$ fixes at least one point, then it is contained in $P\rtimes C_{q-1}$ for some $P\in\Syl_3(G)$. From the structure of this point stabiliser then follow the first three bullet points and the second last. This can also be verified with the description given in \Cref{appendixA}.  Hence, in all other case there are no fixed points. The free actions follow from the fact that, in the cases under consideration, each divisor $d$ of $k$  is again relatively prime with $q^3(q-1)$, and hence, $x^d$ has no fixed points, too. Finally, the fact that no other orders show up is due to the structure of the maximal subgroups, in particular the normalisers of the cyclic Hall subgroups as given above. 
\end{proof}
\color{black}

\begin{lem}\label{lem:noncentralorder3}
Let $\alpha,\beta\in \Omega$ be distinct. Let $P_\alpha$ and $P_\beta$ be the corresponding Sylow $3$-subgroups. Let $x\in [P_\alpha,P_\alpha]\setminus Z(P_\alpha)$ and $y\in [P_\beta,P_\beta]\setminus Z(P_\beta)$. Suppose $\langle x,y \rangle$ does not stabilise $B(\alpha,\beta)$. Then $x$ does not stabilise both blocks $B(\alpha,y\cdot\alpha)$ and $B(\alpha,y^{-1}\cdot\alpha)$ simultaneously.  
\end{lem}

\begin{proof}
We again use the notation of \Cref{appendixA} and put $\alpha=\infty$ and $\beta=O$. Clearly, the statement is trivial if $x$ does not stabilise $B(\infty,O)$ and $y$ does. 
Hence, we may assume that $y$ does not stabilise $B(\infty,O)$. So, up to conjugation, we may assume that $y=(0,b',1)_O$, and $x=(0,a',a'')_\infty$. Then \[y\cdot \infty=\begin{pmatrix} 
 1 & 0 & 0 & 0 & 0 & 0 & 0\\
 b'^\theta & 1 & 0 & 0 & 0 & 0 & 1\\
 -1+b' & 0 & 1 & -b' & 0 & -1 & b'^\theta\\
  -1 & 0 & 0 & 1 & 0 & 0 & 0 \\
  1+b'^{\theta+1} & b' & 0 & 1 &1 & -b'^\theta & 1+b' \\
 -b' & 0 & 0 & 0 & 0 & 1 & 0 \\
0 & 0 & 0 & 0 & 0 & 0 & 1 
 \end{pmatrix}\begin{pmatrix} 1\\ 0\\0\\0\\0\\0\\0\end{pmatrix}=\begin{pmatrix} 1\\b'^\theta\\-1+b'\\-1\\1+b'^{\theta+1}\\-b'\\0\end{pmatrix}\]\[=p\left(\frac{1-b'}{1+b'^{\theta+1}},\frac{-b'^\theta}{1+b'^{\theta+1}},\frac1{1+b'^{\theta+1}}\right).\]
 Now, $y^{-1}=(0,-b',-1)$ and so we find that \[y^{-1}\cdot\infty=p\left(-\frac{1+b'}{1+b'^{\theta+1}},\frac{b'^\theta}{1+b'^{\theta+1}},\frac1{1+b'^{\theta+1}}\right).\]
On the other hand, one easily calculates that $x$ maps the block $B_{c,c''}$ to the block $B_{c,c''+a''-a'c}$, and hence, the block stabilised by $x$ is  $B_{a''a'^{-1},c''}$. Clearly at most one of the values $\frac{1-b'}{1+b'^{\theta+1}}$ or $-\frac{1+b'}{1+b'^{\theta+1}}$ is equal to $a''a'^{-1}$, which proves the lemma.  
\end{proof}

\begin{lem}\label{lem:noncentralorder3bis}
Let $\alpha,\beta\in \Omega$ be distinct. Let $P_\alpha$ and $P_\beta$ be the corresponding Sylow $3$-subgroups. Let $x\in [P_\alpha,P_\alpha]\setminus Z(P_\alpha)$ and $y\in [P_\beta,P_\beta]\setminus Z(P_\beta)$. Suppose both $x$ and $y$ stabilise $B(\alpha,\beta)$. Then each block distinct from $B(\alpha,\beta)$ and fixed by $x$ intersects a unique block distinct from $B(\alpha,\beta)$ and fixed by $y$.
\end{lem}

\begin{proof} 
We use the notation of \Cref{appendixA}. Clearly, the elements $(0,a,0)_\infty$, which are precisely those elements of $P_\infty$ stabilising $B(\infty,O)$, and which, apart from the identity, all belong to $[P_\infty,P_\infty]\setminus Z(P)_\infty)$, stabilise the blocks $B_{(0,a'')}$, containing points of the form $p(0,a',a'')$. Likewise,  $(0,a,0)_O$ stabilises all lines containing points of the form $q(0,b',b'')$. \Cref{calculation} implies that the common points of the sets \[\{p(0,a',a'')\mid a'\in\mathbb{F}_q, a''\in\mathbb{F}_q^\times\}\mbox{ and }\{q(0,b',b'')\mid b'\in\mathbb{F}_q, b''\in\mathbb{F}_q^\times\}\] are the points $p(0,a''^{\theta-1},a'')$. It follows that each block $B_{(0,a'')}$, $a''\neq 0$, contains a unique such point. This proves the assertion.
\end{proof}

A \emph{central element $x$ of order $3$} in $G$ is by definition one that is contained in $Z(P)$, for some Sylow $3$-subgroup $P$ of $G$.  
\begin{lem}\label{lem:centralelements}
A central element of order $3$ of $G$ does not stabilise any block, and is not contained in the normaliser of any (cyclic) Hall subgroup. 
\end{lem}
\begin{proof}
The first assertion follows from \Cref{actionofP}. From the structure of $N_G(M_i)$, $i\in\{-1,0,1\}$ one sees that each element $x$ of order $3$ of $N_G(M_i)$ is contained in a cyclic subgroup of order $6$. Hence, there must exist an involution $\eta$ such that $x\eta$ has order $6$. Clearly $\eta$ has to fix the fixpoint of $x$ which by Lemma \ref{lem:order-fix} is unique. Hence, with the notation of \Cref{appendixA}, we can write $\eta=h(-1)$ and $x=(0,0,b'')_\infty$, for some $b''\in\mathbb{F}_q$. But then one calculates that $x\eta$ maps $p(a,a',a'')$ to $p(-a,a',-a''-b'')$, and hence, $x\eta$ has order $2$, a contradiction. 
\end{proof}

We now show that the maximal subgroups $N(M_0)$ are the counterparts of  the stabilisers of (real) triangles in the unitary groups $\U_3(q)$.
\begin{lem}\label{Hall3blocks}
The maximal subgroup $N(M_0)$ of $G$, for a Hall subgroup $M_0$ of order $\frac12(q+1)$, stabilises three mutually disjoint blocks, each of them stabilised by every member of $N(M_0)$ with an order not divisible by $3$, and permuted around in a cycle by each subgroup of order $3$ of $N(M_0)$. 
\end{lem}

\begin{proof}
The subgroup $N(M_0)$ is the normaliser of a Klein four group $V=:\{1,\eta_1,\eta_2,\eta_3\}$. As such, each member of it permutes the three fix blocks $B_i$ of the $\eta_i$, respectively, $i=1,2,3$, and $N(M_0)$ does not contain further involutions. Each $\eta_i$ itself, $i=1,2,3$, centralises $V$ and hence stabilises each block $B_i$, $i=1,2,3$. It follows that $B_1\cap B_2=\varnothing$ as otherwise $\eta_1$ fixes a unique point of $B_2$, leading to $2|q$, a contradiction. Let $g\in N(M_0)$ be an arbitrary member that does not stabilise all of $B_1,B_2,B_3$.      Suppose first that $g$ stabilises $B_1$ and interchanges $B_2$ and $B_3$. Then $|g|$ is even. By the structure of $N(M_0)$, we see that $|g|=2k$, for $k$ odd. Hence,  $g^k=\eta_1$ and $g^k$ stabilises both $B_1,B_2$, a contradiction. Now assume that $g$ acts as an order $3$ permutation on $\{B_1,B_2,B_3\}$. Then clearly $|g|\in\{3,6\}$ by the structure of $N(M_0)$ again (alternatively, note that any element whose order divides $3$ is contained in the centraliser of an order $3$ element, which is on its turn contained in a point stabiliser, from which the claim also follows). If $|g|=3$, then $g$ has a unique fixed point $\alpha$ by Lemma \ref{lem:order-fix}, which we claim is not contained in $B_1\cup B_2\cup B_3$. Indeed, suppose for a contradiction that $\alpha\in B_1$. Then by the disjointness of the blocks $B_1,B_2,B_3$, the block $B_1$ is stabilised by $g$. But we proved above that in this case also $B_2$ and $B_3$ have to be stabilised, contradicting \Cref{actionofP} and the fact that $B_2$ and $B_3$ are disjoint from $B_1$. Hence, the claim follows and \Cref{actionofP} implies that $g$ cycles $B_1,B_2,B_3$ around. If $|g|=6$, then $g^2$ cycles these blocks around, and hence also $g$.  
\end{proof}

\subsection{Structural subgroups}

Structural subgroups are subgroups of the maximal subgroups $(i)$ $N_G(P)$, for $P \in \mathsf{Syl}_3(G)$, and $(ii)$ $C_G(\eta)$, where $\eta$ is an involution. 

\emph{$(i)$}: The normalisers of Sylow $3$-subgroups each stabilise a point. Sylow $3$-subgroups consist of unipotent elements, have order $q^3$, their center is elementary abelian of order $q$, they have class $3$ and the Frattini subgroups equal the commutator subgroups and contain the center.

The stabiliser of two points $\alpha,\beta$ is a cyclic group $W$ of order $q-1$ and $N_G(P)=P\rtimes W$, where $P\in\Syl_3(G)$ and fixes the point $\alpha$. The group $W$ acts by conjugation on the set of elements of $[P,P]$ that stabilise $B(\alpha,\beta)$ in exactly two orbits, mapped onto each other by taking inverses. 

\emph{$(ii)$}: The centraliser of an involution $r$ is isomorphic to $\langle r\rangle \times \PSL_2(q)$ and these are the block stabilisers.

The following table lists the maximal subgroups of $G$ up to conjugacy and their geometric properties.

\begin{center}
\renewcommand{\theadfont}{\normalsize\bfseries}
\begin{tabular}{ r | l | l}
  \thead{Maximal subgroup} & \thead{Geometric prop.} & \thead{Orders}\\
\hline
 $N_G(P)$, $P \in \mathsf{Syl}_3(G)$ & Point stabiliser & $2$, $3$, $6$, $9$\\
  $C_G(r)$, $r$ involution & Block stabiliser & $2$, $3$, $6$, divisors of $q-1$, $(q+1)/2$\\ 
\hline
 $^2\mathsf{G}_2(q_0), \: q_0^p = q$, & Subunital stabiliser  & \\
  $p$ prime & &  \\
\hline
 $N_G(M_{\pm 1})$ & & $2$, $3$, $6$, divisors of $q \pm \sqrt{3q}+1$\\
   $N_G(M_0)$ & Triangle stabiliser & $2$, $3$, $6$, divisors of $(q+1)/2$
\end{tabular}
\end{center}

Subgroups of the first two rows will be referred to as \emph{structural subgroups}. Subgroups of maximal subgroups of the form $^2\mathsf{G}_2(q_0)$ are called \emph{subunital stabilisers} or \emph{subfield subgroups}. All other subgroups will be referred to as \emph{small subgroups}.
\subsection{Classification of elements in $^2\mathsf{G}_2(q)$ and their normalisers}
\label{subsection: ree structure}

A list of conjugacy classes in $G$ is given by the following table, where the left column denotes either the order of an element, or a multiple $m$ of that order, and the right column denotes either the number of conjugacy classes of elements of that given order, or the total number of conjugacy classes of elements whose order divides $m$, respectively.
\begin{table}[h]
\begin{center}
\begin{tabular}{ c c c  }
  & Order & Number of conjugacy classes\\
 \cmidrule{2-3}
\hspace{-1em}\ldelim\{{4}{*}[Order divides] & $q+\sqrt{3q}+1$ & $\frac{q+\sqrt{3q}}{6}$  \\ 
 & $q-\sqrt{3q}+1$ & $\frac{q-\sqrt{3q}}{6}$ \\
 & $\frac{q+1}{2}$ & $\frac{q-3}{6}$\\
 & $\frac{q-1}{2}$ &  $\frac{q-3}{2}$\\
\cmidrule{2-3}
\hspace{1.3em}\ldelim\{{5}{*}[Order is]  & $9$ & $3$\\
 & $6$ & $2$\\
 & $3$ & $3$\\
& $2$ & $1$\\
& $1$ & $1$ 
\end{tabular}
\end{center}
\caption{Number of conjugacy classes in $^2\mathsf{G}_2(q)$\label{table1}}
\end{table}

Given this table, the facts listed in \Cref{sec:generalities} and the table of maximal subgroups of $G$, the following proposition is obvious. 

\begin{prop}\label{centralisers}
Let $x$ be a nontrivial element of $G$. Then $C_G(x)$ is given by the following subgroups.
 \begin{eqnarray*}
C_{G}(x) =
\begin{cases}
C_{{q}-1} & |x| \mid \frac{q-1}{2} \\
C_2^2 \times C_{\frac{q+1}{4}} & |x| \mid \frac{q+1}{4} \\
C_{q\pm\sqrt{3q}+1} & x \in M_i, i = \pm 1\\
P & x \in Z(P), P \in \Syl_3(G)\\
[P,P] \times \left<r \right>, \: |r| = 2 & x\in [P,P] \setminus Z(P), P \in \Syl_3(G)\\
Z(P) \cdot \left<x\right> & x \in P \setminus [P,P], P \in \Syl_3(G)\\
\left<x\right> \times \PSL_2(q) & |x|=2
\end{cases}
\end{eqnarray*}
\end{prop}

\begin{lem}\label{lem:conjugate in smaller ree group}
Suppose that $x,y \in G$, where $x$ is not an involution, and that $\left<x, x^y\right>$ is not a structural subgroup and is not contained in a subgroup isomorphic to $^2\mathsf{G}_2(3)'\cong\PSL_2(8)$. Then $x$ and $x^y$ are conjugate in $\left<x,x^y\right>$.
\end{lem}
\begin{proof}
This follows readily from the description of the conjugacy classes. For the Hall subgroup stabilisers, this is a consequence of the fact that they are normalised by a cyclic group of order $6$, whereas the corresponding denominator in the right column of \Cref{table1} is exactly $6$. For subfield subgroups $H$, this follows from the observation that the number of conjugacy classes of elements of $H$ of a given order in $H$ is equal to the number of such conjugacy classes in $G$, together with the fact that, if this number is strictly larger than $1$, which happens only if $|x|\in \{3,6,9\}$, then the following holds. If $|x|=3$, then one class is given by the centre $Z(P')$ of a Sylow 3-subgroup $P'$ of $H$, where $Z(P')\leq Z(P)$, for some Sylow $3$-subgroup $P$ of $G$, and the two other classes are related by inverse (so all three classes are present in every subfield subgroup). If $|x|=6$, then the two classes are related by inverse. For $|x|=9$, two of the classes are related by inverse and the third class has, with the notation of \Cref{appendixA}, representative $(1,-1,0)_\infty$ and hence exists in every Ree subgroup, except in the smallest derived case $^2\mathsf{G}_2(3)'\cong\PSL_2(8)$. 
\end{proof}

\begin{lem}\label{lem:centralizer ree subgroup inequality}
Let $x$ be an element which is not an involution in $G$, and suppose that $x \in R_1 \cap R_2$ where $R_1$ and $R_2$ are Ree subgroups of $G$, distinct from $^2\mathsf{G}_2(3)'$. If $|C_{R_1}(x)| \leq |C_{R_2}(x)|$, then $|R_1| \leq |R_2|$.
\end{lem}

\begin{proof}
This follows from \Cref{centralisers} by noting that the orders of the centralisers depend proportionally on the order of the Ree subgroup. Also, elements of a given Ree subgroup $R_1$ stay ``in the same kind of class'' when considered in another Ree subgroup $R_2$; that is, if $C_{R_1}(x)$ for some member $x\in R_1$ is read on the $n$th line of the displayed cases in \Cref{centralisers}, then on the same line one reads $C_{R_2}(x)$. This follows from the fact that, if $R_i$ is defined over the field $\mathsf{F}_{q_i}$, the order $|x|$ of any member $x\in R_1\cap R_2$ divides $q_i+1$ if, and only if, it divides $q+1$; it divides $q_i-1$ if, and only if, it divides $q-1$, and it divides $q_i^2-q_i+1$ if, and only if, it divides $q^2-q+1$. (These follow from easy divisibility conditions possibly using elementary cyclotomic polynomial theory.)
\end{proof}

\subsection{Basic group theory lemmas}

We conclude this section with two basic group theory lemmas that will be useful throughout this article.

\begin{lem}\label{lem:cyclic subgroup reduction}
Let $K$ be a finite cyclic group generated by $\{x_1,x_2, \ldots, x_k\}$. Then  there exists natural numbers $m_1, \ldots, m_{k-1}$ so that 
$$
K = \left<x_1^{m_1} \ldots x_{k-1}^{m_{k-1} }x_k \right>.
$$
\end{lem}
\begin{proof}
This proof follows from \cite[Lemma 4.1]{G08} and induction.
\end{proof}

The proof of the following lemma is immediate.
\begin{lem}\label{lem:semi-directproductprop}
Let $K=C_n\rtimes C_m$ be a semi-direct product of two cyclic groups with $(n,m)=1$. If the prime $p$ divides $m$, then the elements $g\in K$ with $g^{pn}=1$  form the subgroup $C_n\rtimes C_p$. If $p=3$, then for any pair $g,h\in C_n\rtimes C_3$ of elements of order $3$ we either have $|gh|=3$ and $(gh)^n=1$, or $|gh^{-1}|=3$ and $(gh^{-1})^n=1$. 
\end{lem}

We now embark on the proof of Theorem \ref{main_thm1}. We show that every vertex of $\widetilde{X}_5(G)$, for $G={^2\mathsf{G}}_2(q)$, is connected to a redundant vector. We prove as much as we can for shorter vectors, in particular for vectors of length 3. Globally we show that we can transform any vector $(x,y,z,v,w)\in V_5(G)$ with a finite number of Nielsen moves to a redundant vector. Obviously we may assume that $(x,y,z,v,w)$ is not redundant itself, in particular, none of $x,y,z,v,w$ is the identity.

\section{Reduction away from involutions and unipotent elements}\label{sec:reductions}
This first lemma says that we may connect any generating $3$-tuple $(x,y,z)$ to a generating $3$-tuple $(x,y',z')$ such that $y'$ and $z'$ are not involutions and where $y'$ and $z'$ do not normalise the cyclic subgroup $\left<x\right>$.

\begin{lem}\label{lem:normalizer_involution}
Let $G$ be a finite simple group, and let $(x,y,z) \in V_3(G)$. Then $(x,y,z)$ is connected to $(x, y', z')$, where $y',z' \notin N_G(\left<x\right>)$. Moreover, we may assume that $y'$ and $z'$ are not of order $2$.
\end{lem}
\begin{proof}
\emph{Step 1:} We observe that $y$ and $z$ cannot both be in the normaliser $N_G(\left<x\right>)$, since $G$ does not normalise $\left<x\right>$. We note that if $y \in N_G(\left<x\right>)$ and $z \notin N_G(\left<x\right>)$, then $yz \notin N_G(\left<x\right>)$. Thus, we can connect $(x,y,z) \to (x,yz,z)$, where $yz,z \notin N_G(\left<x\right>)$. If $yz$ is an involution, we proceed with Step 2 with $yz$ in place of $y$.

\emph{Step 2:} If $y \notin N_G(\left<x\right>)$, but $y$ is an involution, we claim that $xy$ is not an involution. If it were, then $xyxy = 1$, which implies $x^{-1} = yxy = x^y$. Hence, $x^y \in \left<x\right>$. That implies $y$ normalises $\left<x\right>$, which is a contradiction. Thus, if $y$ has order $2$, we may connect $(x,y,z)$ to $(x,xy,z)$ where $xy$ does not have order $2$, and since $y \notin N_G(\left<x\right>)$, it follows that $xy \notin N_G(\left<x\right>)$. We can do the same step with $z$ instead of $y$ and arrive at the same assertion. 
\end{proof}

Now we show that we can connect every quintuple to a quintuple containing at least one element that is not \emph{unipotent}, that is, which does not fix a unique point (such elements are characterised by the property that they have order $3,6$ or $9$). 
Consider the following condition on a quintuple $(x,y,z,v,w)$.

\begin{itemize}
\item[\textbf{(*)}] Either at least one of the elements $x,y,z,v,w$ has order distinct from $2,3,6$ and $9$, or none of the elements are involutions and the subgroups $\<x,y,z,v\>$, $\<x,y,z,w\>$, $\<x,y,v,w\>$ and $\<x,z,v,w\>$ are subfield subgroups over subfields of size at least $27$, and the subgroups $\<x,x^y\>$, $\<x,x^z\>$, $\<x,x^v\>$ and $\<x,x^w\>$ are not structural subgroups. 
\end{itemize}

We first show that, in order to prove that for a given generating quintuple $(x,y,z,v,w)$, we may assume that the subgroup $\<x,y,z,u\>$ is a subfield subgroup over a subfield of size at least 27, it suffices to show that it is not structural. 

\begin{lem}\label{smallthenredundant}
Let $(x,y,z,v,w)\in V_5(G)$. Then, if $\<x,y,z,v\>$ is small, then $(x,y,z,v,w)$ is connected to a redundant $5$-tuple.
\end{lem}

\begin {proof}
Small subgroups of $G$ lie, by their very definition, in $N_G(M_i)$, which is the normaliser of the cyclic Hall subgroup $M_i$, $i= -1,0,1$, of order $q-\sqrt{3q}+1$, $ \frac{q+1}{4}$, and $q+\sqrt{3q} + 1$, respectively.  For $i\in\{-1,1\}$ all subgroups of $N_G(M_i)$ are either cyclic or $2$-generated since then $N_G(M_i)$ is cyclic-by-cyclic. The structure of $N_G(M_0)$ is given by $(V \rtimes M_0) \rtimes \left<t\right>$ where $V$ is an elementary abelian group of order $4$ and $t$ has order $6$. We note that $t$ acts on $V$ by cyclically interchanging the three involutions. Hence, all subgroups of $N_G(M_0)$ are at most $3$-generated. In all cases $\<x,y,z,v\>$ admits a generating set of $3$ elements. Now, work of Dunwoody \cite{DW} implies that $\left<x,y,z,v\right>$ is connected to a redundant generating $4$-tuple. Hence, $(x,y,z,v,w)$ is connected to a redundant $5$-tuple.\end{proof}

\begin{lem}\label{lem:eliminating small subgroups}
Let $(x,y,z,v,w)$ be a generating $5$-tuple such that $\<x,y,z,u\>$ is either small or isomorphic to either $^2\mathsf{G}_2(3)'$ or $^2\mathsf{G}_2(3)$. Then $(x,y,z,v,w)$ is connected to a redundant $5$-tuple
\end{lem}
\begin{proof}
For $\<x,y,z,u\>$ small, this follows directly from \Cref{smallthenredundant}. Observe that ${^2\mathsf{G}}_2(3)\cong \mathsf{P\Gamma L}_2(8)\cong\PSL_2(8)\rtimes 3$.  
 Hence, $H=:\<x,y,z,u\>\in\{\PSL_2(8),\PSL_2(8)\rtimes 3\}$. If $H\cong \PSL_2(8)$, then by \cite{E93}, $(y,z,v,w)$ is connected to a redundant $4$-tuple, a contradiction. So we may assume $H\cong\PSL_2(8)\rtimes 3$ and we set $H=H_0\rtimes\<\xi\>$, with $|\xi|=3$. We may assume $w\notin H_0$. Then for some $i\in\{0,1,2\}$, $w^iy\in H_0$. This yields Nielsen moves to connect $(y,z,v,w)$ to $(y',z',v',w)$ with $y',z',v'\in H_0$. By the classification of subgroups of $\PSL_2(8)$, we either have $\<y',z',v'\>\cong\PSL_2(8)$, in which case $(y',z',v')$ is connected to a redundant triple by \cite{E93} again, or $\<y',z',v'\>$ is solvable and generated by at most two elements, implying by \cite{DW} that we again can connect $(y',z',v')$ to a redundant triple.
\end{proof}

\begin{prop}\label{condition(*)}
Every quintuple of generating elements of $G$ containing only elements of orders $2,3,6,9$ is connected in $\widetilde{X}_5(G)$ to a quintuple satisfying\emph{ Condition (*)}. 
\end{prop}

\begin{proof} Let the quintuple be given by $(x,y,z,v,w)$. We first note that not all of $x,y,z,v,w$ are involutions, as, if they mutually commuted, they would generate an elementary abelian subgroup of order $32$, a contradiction. Hence, if $x$ and $y$ are non-commuting involutions, the quintuple $(xy,y,z,v,w)$ contains an element $xy$ of order at least $3$. So, without loss of generality, we may assume that $x$ is not an involution. 

By \Cref{lem:normalizer_involution}, we may assume that all of $y,$ $z$, $v$ and $w$ are not involutions. Hence, each of $x,y,z,v,w$ fixes a unique point, which we denote by $\alpha, \beta, \gamma, \delta, \xi$, respectively. Evidently, we do not have $\alpha=\beta=\gamma = \delta=\xi$. If $\alpha=\beta$, then we may assume without loss of generality that $\gamma\neq\alpha$ and we replace $(x,y,z,v,w)$ with $(x,y^z,z,v,w)$. Repeating this if necessary, 
this produces a quintuple with $\alpha\notin\{\beta,\gamma,\delta,\xi\}$.  This already implies that we may assume that neither $\<x,y\>$, $\<x,z\>$, $\<x,v\>$ nor $\<x,w\>$ are contained in a point stabiliser. We divide into cases. 

\textbf{Case 1.} \emph{Suppose at least one of $x,y,z,v,w$ has order $9$.} We may assume that $x$ has order $9$. Since elements of order $9$ only appear in subfield subgroups and point stabilisers, $\<x,y\>$ is a subfield subgroup. Since $x^y$ only fixes the point $y\cdot\alpha\neq\alpha$, the group $\<x,x^y\>$ is a subfield subgroup. Similarly all of $\<x,z\>$, $\<x,x^z\>$, $\<x,v\>$, $\<x,x^v\>$, $\<x,w\>$ and $\<x,x^w\>$ are subfield subgroups and the assertion follows from the fact that overgoups of subfield subgroups are subfield subgroups, as is apparent from the list of maximal subgroups. . 

\textbf{Case 2.} \emph{None of $x,y,z,v,w$ has order $9$, but at least one of them is central.}  We may suppose that $x$ is central. By \Cref{lem:centralelements}, the subgroups $\<x,y\>$, $\<x,x^y\>$, $\<x,z\>$, $\<x,x^z\>$, $\<x,v\>$, and $\<x,x^v\>$ are subfield subgroups. The assertion again follows.

\textbf{Case 3.} \emph{None of $x,y,z,v,w$ have order $9$ or are central, and one of $x,y,z,v,w$ has order $6$.}  We may assume $|x|=6$. Note that $x$ stabilises a unique block, which contains $\alpha$. Suppose both $x$ and $y$ stabilise $B(\alpha,\beta)$. Then, since we may then assume without loss of generality that $z$ does not stabilise $B(\alpha,\beta)$, one of $z$ or $z^{-1}$ maps $\beta$ outside $B(\alpha,\beta)$. Suppose without loss of generality that $z\cdot\beta\notin B(\alpha,\beta)$. Then we replace $(x,y,z,v,w)$ with $(x,y^z,z,v,w)$. Hence, we may assume not both $x$ and $y$ stabilise $B(\alpha,\beta)$.  Similarly not both $x$ and $z$ stabilise $B(\alpha,\gamma)$, not both $x$ and $v$ stabilise $B(\alpha, \delta)$, and not both $x$ and $w$ stabilise $B(\alpha,\xi)$. Then $\<x,y\>$, $\<x,z\>$, $\<x,v\>$ and $\<x,w\>$ are not structural. Hence, $\<x,y,z,v\>$ is not structural. But it is not small either by \Cref{smallthenredundant}; hence, it is a subfield subgroup. Similarly, the subgroups $\<x,y,z,w\>$, $\<x,y,v,w\>$ and $\<x,z,v,w\>$ are subfield subgroups. 

Suppose $x$ and $x^y$ stabilise the same block $B$. Then $B$ contains $\alpha$ and $y\cdot\alpha$. Then the block $B(\alpha,y^{-1}\cdot\alpha)$ is distinct from $B$ and hence not stabilised by $x$. It follows that $\<x,{^yx}\>$  is not structural. Replacing $y$ with $y^{-1}$, which is allowed by the Nielsen moves, we obtain that $\<x,x^y\>$ is not structural. On the other hand, if $x$ and $x^y$ do not stabilise the same block, then $\<x,x^y\>$ was not structural in the first place. Similarly, either $\<x,x^z\>$ or $\<x,{^zx}\>$, either $\<x,x^u\>$ or $\<x,{^vx}\>$, and either $\<x,x^w\>$ or $\<x,{^wx}\>$ is not structural. Then an analogous argument as with $\<x,x^y\>$ yields (*).

\textbf{Case 4.} \emph{All of $x,y,z,v,w$ have order $3$ and are not central.} Here, each of $x,y,z,v,w$ stabilises $q$ blocks through their respective fixed point. We claim that we may assume that $x$ and $y$ do not share a common fixed block. Indeed, assume all of $x,y,z,u,w$ pairwise share common fixed blocks. Evidently, this cannot be the same block for every pair; hence, we may assume that $\gamma$ is not contained in $B(\alpha,\beta)$. \Cref{lem:noncentralorder3bis} then implies that $\beta$ is the unique point on $B(\beta,\gamma)$ whose  join to $\alpha$ is fixed by $x$. Consequently, $y^z$ does not fix any block fixed by $x$, and since we may replace $y$ with $y^z$, the claim follows. 

Hence, $\<x,y\>$ is not structural and, by \Cref{lem:noncentralorder3}, and replacing $y$ by its inverse if necessary, the subgroup $\<x,x^y\>$ is not structural. We now claim that we may assume that $\<x,x^z\>$ is not structural. Indeed, this follows from \Cref{lem:noncentralorder3} if $x$ does not fix the block $B(\alpha,\gamma)$. So assume that $x$ stabilises $B(\alpha,\gamma)$. Since $\<x,y\>$ is not structural, the subgroup $H:=\<x,y,u,w\>$ is not structural and so, by \Cref{lem:eliminating small subgroups}, $H$ is isomorphic to $\mathsf{^2G_2}(q_0)$ for some $3$-power $q_0\geq 27$. By field restriction it follows that $H$ acts doubly transitive on a subunital $\Omega_0$ of $\Omega$, which contains $\alpha$ as $\alpha$ is the unique fixed point of $x$. Let $u\in H$ be an element of order $9$ fixing $\alpha$. Then, replacing $z$ by $z'=z^u$, as is allowed with Nielsen moves, \Cref{actionofP'} implies that $x$ does not fix $B(\alpha,u\cdot\gamma)$. So, again by \Cref{lem:noncentralorder3}, and replacing $z'$ by its inverse if necessary, $\<x,x^{z'}\>$ is not structural. The claim follows. Similarly we may assume that $\<x,x^u\>$ and $\<x,x^w\>$ are not structural. But then non of $\<x,y,z,v\>$,  $\<x,y,z,w\>$, $\<x,y,v,w\>$ and $\<x,z,v,w\>$ is structural, and hence, by \Cref{smallthenredundant}, they are all subfield subgroups. \color{black}

The proof is complete.\end{proof}

\begin{prop}\label{unipotentconnectsto2or0}
Every quintuple $(x,y,z,v,w)$ of generating elements of $G$ containing only elements of orders $2,3,6,9$ is connected in $\widetilde{X}_5(G)$ to either a redundant quintuple, or to a quintuple containing at least one element of order distinct from $2,3,6$ or $9$.
\end{prop}

\begin{proof}
By Proposition~\ref{condition(*)}, we may assume that  $H:=\<x,y,z,u\>$ is a subfield subgroup $\mathsf{^2G_2}(q_0)$, with $q_0\geq 27$.

We now claim that the coset $wH$ contains at least one element of order distinct from $2,3,6$ or $9$. We proceed by distinguishing the various possible orders of $w$. We perform explicit computations using the matrices given in \cite[Section~9.2.4]{TTM07}, see \Cref{appendixA} below.  Those $7\times 7$ matrices over $\mathbb{F}_q$ describe the action of the full Sylow 3-subgroups $P_\infty$ and $P_O$, which fix the two points $\infty$ and $O$, respectively. Recall that each Sylow 3-subgroup is described by triples $(a,a',a'')_\infty$ and $(a,a',a'')_O$, respectively, with $a,a',a''\in\mathbb{F}_q$, Also, $Z(P_\infty)=\{(0,0,a'')_\infty\mid a''\in\mathbb{F}_q\}$, that $[P_\infty,P_\infty]=(0,a',a'')_\infty\mid a',a''\in\mathbb{F}_q\}$ (and similarly for $P_O$), and  the stabiliser of the block through $\infty$ and $O$ in $P_\infty$ and $P_O$ is $\{(0,a',0)_\infty\mid a'\in\mathbb{F}_q\}$ and $\{(0,a',0)_O\mid a'\in\mathbb{F}_q\}$, respectively. The $2$-point stabiliser $G_{\infty,O}$  acts by conjugation on $P_\infty$ as $(a,a',a'')_\infty\mapsto (at^\theta,a't^{\theta+1},a''t^{\theta+2})_\infty$. Recalling the definition of trace of an element of $G$ from \Cref{appendixA}, an element of order $3$ or $9$ has trace 1 and an element of order $2$ or $6$ has trace $-1$ by \Cref{tracelemma}.

We may assume that $w$ fixes $\infty$. Let $U$ be the group of central elements of an arbitrary Sylow $3$-subgroup of $H$ not fixing $\infty$. Then, by the $2$-transitivity of $(G,\Omega)$, we may assume that $U$ fixes $O$.
A generic member $u$ of $U$ is $(0,0,d)_O$, $d\in\mathbb{F}_q$, but note that not all elements of $\mathbb{F}_q$ provide elements of $U$; only $q_0$ do. However, we may assume that for those $q_0$ values, the traces of $wu$ and $wu^2$ are in $\{1,-1\}$. We consider all possible orders for $w$ and prove that under this assumption, each lead to a contradiction.

\begin{itemize}
\item \emph{Assume that $w$ has order $9$.} 
By the action of the $2$-point stabiliser $G_{\infty,O}$, we may assume that $w=(1,b,c)_\infty$, for some $b,c\in\mathbb{F}_q$. One can now explicitly calculate the traces of the matrices corresponding to the products $wu^{\pm1}$ and obtain
\[1+d^2(c^2-c^\theta+b^\theta+b^{\theta+1}-b^2-b-1)\mp [d+d^\theta(b^2-b-c-c^\theta)],\]
where the minus sign corresponds to $wu$ and the plus sign to $wu^2=wu^{-1}$. Since these traces belong to $\{1,-1\}$, the difference of the traces of $wu$ and $wu^{-1}$ belongs to $\{0,1,-1\}$. Hence, we obtain,
setting $b^2-b-c-c^\theta=f$ (possibly $f=0$), that $d+fd^\theta=\varepsilon$, with $\epsilon\in\{0,1,-1\}$. This implies $d^\theta+f^\theta d^3=\varepsilon^\theta$, yielding \[(f\varepsilon^\theta-\varepsilon)+d-f^{\theta+1}d^3=0.\] The latter can only be true for at most three nonzero values of $d$, whereas we have $q_0$ values of $d$, a contradiction. 

\item  \emph{Assume that $w$ has order $3$ and does not fix every block joining its fixed point $\infty$ to a point of the subunital stabilised by $H$.} Similarly as in the previous case, we may assume that $w=(0,b,1)_\infty$, and both $wu$ and $wu^2$ have trace either $1$ or $-1$ (with $u$ as above). This implies
\[\begin{cases}
1+d^2 +b^{\theta+1} d^2 +d^\theta\in\{1,-1\},\\
1+d^2 +b^{\theta+1} d^2 -d^\theta\in\{1,-1\},
\end{cases}\]
which obviously implies $d\in\{0,1,-1\}$. 
Since this must be true for all $q_0$ values for $d$, this is again a contradiction.

\item \emph{Assume that $w$ has order $3$ and fixes every block joining its fixed point $\infty$ to a point of the subunital stabilised by $H$.}  Here we can take $w=(0,b,0)_\infty$ and $u=(0,0,d)_O$. Calculating the trace of $wu$, we obtain $1+b^{\theta+1}d^2$, which is only equal to $1$ if one of $b$ or $d$ is $0$, a contradiction. If the trace of $wu$ is $-1$ then $b^{\theta+1}d^2=1$. Since there are at most two solutions for $d$, we again reached a contradiction.
\item \emph{Assume finally that $w$ has order $6$.}  Any order $6$ element fixing $\infty$ can be obtained by multiplying $(0,b,0)_\infty$ with the involution $\eta$ that fixes the block $B(\infty,O)$ pointwise, and which  is given by the diagonal matrix $\mathsf{Diag}({1, 1, - 1, - 1, 1, 1, - 1})$. The set of elements of the form $(a,0,c)_\infty$ acts transitively on the set of blocks through $\infty$, and conjugating with an appropriate element of the $2$-point stabiliser $G_{\infty,O}$, we may assume $a=1$. Then one calculates that the trace of 
\[(1,0,c)_\infty^{-1}(0,b,0)_\infty\eta(1,0,c)_\infty(0,0,\pm d)_O\] 
is equal to 
\[1 \mp bd +d^2-b^2d^2+b^{\theta+1}d^2+cd^2+c^2d^2-c^\theta d^2 \pm b^\theta d^\theta \pm bc d^\theta\]
where the plus and minus signs should be read consistently. Hence the difference of these traces is equal to 
\[bd-b^\theta d^\theta-bcd^\theta\] 
Since both traces are either $1$ or $-1$ we obtain that $bd-b^\theta d^\theta-bcd^\theta = k$ where $k\in \{0,1,-1\}$.
Hence also $b^\theta d^\theta-b^3 d^3-bcd^3 = k$. Since $b\neq 0$ we can compute $d^\theta$ from the second equation and substitute in the first to obtain either a linear or a cubic equation in $d$ which can be satisfied for at most $3<q_0$ values of $d$, again a contradiction.
 \end{itemize}
Hence, the claim is proved. Let $h\in H$ be such that the order of $wh$ is not $2,3,6$ or $9$.  Then we can do Nielsen moves to transform $(w,x,y,z)$ into $(wh,x,y,z)$ and the proposition follows.
\end{proof}

\section{Eliminating the structural subgroups}\label{sec:structural}
In this section, we show that we can connect a generating $5$-tuple $(x,y,z,v,w)$, where $|x| \notin \{2,3,6,9\}$, to a generating $5$-tuple where $\left<x,x^y\right>$,  $\left<x,x^z\right>$, $\<x,x^v\>$ and $\<x,x^w\>$ are not structural subgroups. See  \Cref{subsection: ree structure} for the list of structural subgroups of the Ree groups. We proceed based on the fixed point sets of $x$. 

\subsection{$x$ fixes two points}
This first proposition allows us to assume when given a generating $3$-tuple $(x,y,z)$ where $x$ fixes exactly two points (hence $|x| \neq 2$), that we may connect $(x,y,z)$ to a generating $3$-tuple $(x,y',z')$ where $\left<x,y'\right>, \left<x,z'\right>$ are not structural subgroups. Moreover, we may assume that $|y'|, |z'| \neq 2$. For this proposition, when given a point $\alpha$, we denote $G_\alpha = \text{Stab}_G(\alpha)$, and when given two points $\alpha, \beta$, we write $G_{\alpha, \beta} = G_\alpha \cap G_\beta$. We denote the block through the points $\alpha$ and $\beta$ as $B(\alpha,\beta)$, and write its stabiliser as $G_{B(\alpha,\beta)} = G_B$.

\begin{prop}\label{prop:remove structural subgroups 1}
Let $(x,y,z) \in V_3(G)$ where $x$ fixes two points and $|x| \neq 2$. Suppose that $(x,y,z)$ is not connected to a redundant $3$-tuple. Then $(x,y,z)$ is connected to $(x,y',z')$, where $\left<x,y'\right>$ and $\left<x,z'\right>$ are not structural subgroups. Moreover, neither $y'$ nor $z'$ are involutions.
\end{prop}

\begin{proof}
By \Cref{lem:normalizer_involution}, we may assume that $y,z \notin N_G(\left<x\right>)$ and $|y|,|z| \neq 2$. Suppose that $x$ fixes exactly two points $\alpha$ and $\beta$ and that $H = \left<x,y\right>$ is a structural subgroup. That implies that $H$ is a subgroup of either  $G_\alpha$, $G_\beta$, or $G_{B}$ by Lemma \ref{lem:no-blocks}. We start by assuming $H$ is a subgroup of $G_\alpha$.

\noindent \textbf{Step 1:} If $y$ fixes two points, one of them $\alpha$, we set $y_{1} := y$ and proceed with step 2. Suppose $y$ fixes only the point $\alpha$. In this case, we have $H \leq G_{\alpha}$ where $G_{\alpha} = N_G(P)$, $P \in \Syl_3(G)$. In particular, we have $G_{\alpha} = P \rtimes G_{\alpha, \beta}$, where $G_{\alpha,\beta}$ is cyclic of order $q-1$. 

Let $\sigma$ be the involution fixing the block $G_{\alpha,\beta}$ pointwise and set $Q=P\rtimes\<\sigma\>$. Then $Q$ is the set of elements of order $2$, $3$ or $6$ in $G_{\alpha}$. But $Q$ can also be described as the set of elements of $G_{\alpha}$ not fixing exactly two points. This implies $y\in Q$. 

Hence, if $xy$ were in $Q$, then $xyy^{-1} \in Q$, which would imply that $x \in Q$. This contradicts the fact that $x$ fixes $\alpha$ and $\beta$. We conclude that $xy$ is not in $Q$. In particular, it must have exactly two fix points, one of which is $\alpha$. 
We set $y_{1} := xy$ and proceed with step 2. 

\noindent \textbf{Step 2:} $y_{1}$ fixes two points $\{\alpha, \gamma\}$. \newline
If $\gamma = \beta$, then $H$ is a cyclic subgroup and $H \leq G_{B}$ and we can connect $(x,y_{1}, z)$ to $(x,x^my_{1},z)$, where $\left<x^my_{1}\right>$ generates $H$ (cf.~\Cref{lem:cyclic subgroup reduction}). Hence, $(x,y_{1},z)$ is then connected to a redundant triple, which implies we may assume that $\gamma \neq \beta$.

Since $(x,y_{1},z)$ is a generating $3$-tuple, we have $z \notin G_{\alpha}$.  Hence, $z \cdot \alpha \neq \alpha$. Suppose that $z \cdot \alpha = \beta$ and $z^{-1} \cdot \alpha = \beta$. Then $z$ is an involution (since the square of an element of order $6$ only fixes one point by \Cref{lem:order-fix}), a contradiction. 

Then we define $z^\prime = z$ if $z\cdot\alpha\ne\beta$ and $z'=z^{-1}$ otherwise. Then $z^\prime \cdot \alpha \neq \alpha$ and $z^\prime \cdot \alpha \neq \beta$. Given that $x$ is not an involution, the elements $g_{i} := z^\prime x^i$ for $0 \leq i \leq 2$ are distinct. We see that
\[
y_{1}^{g_i} \cdot( g_i \cdot \alpha) = g_i y_{1} g_i^{-1} g_i \cdot \alpha = g_i y_{1} \cdot \alpha = g_i \cdot \alpha,\]
and similarly, $y_{1}^{g_i} \cdot (g_i \cdot  \gamma) = g_i \cdot \gamma.$ Hence, $y_{1}^{g_i} \in G_{g_i \cdot \alpha, g_i \cdot \gamma}$.  Suppose that $g_i \cdot \gamma = g_j \cdot \gamma$, where $0 \leq i < j \leq 2$. We then have $x^{j-i} \cdot \gamma = \gamma$. Hence, we must have $j=i$. Therefore, there exists an integer  $k \in \{0,1,2\}$ such that $g_k \cdot \gamma \notin \{\alpha, \beta\}$ by the pigeon hole principle. Now we also have $g_k\cdot\alpha=z'\cdot\alpha\notin\{\alpha,\beta\}$.
Hence, $y_{1}^{g_k}$ has two fixed points, neither of which is $\alpha$ or $\beta$. 
If at least one of the two fixed points of $y_1^{g_k}$ is not contained in $B(\alpha,\beta)$, then $\<x,y_1^{g_k}\>$ is not a structural subgroup. Therefore, we may assume that $\<x,y_1^{g_k}\>$ fixes $B(\alpha, \beta)$. 

If at least one of the two fixed points of $y_1^{g_k}$ is not contained in $B(\alpha,\beta)$, then $\<x,y_1^{g_k}\>$ is not a structural subgroup by \Cref{lem:no-blocks}. Therefore, we may assume that $\<x,y_1^{g_k}\>$ stabilises $B(\alpha, \beta)$. In this case clearly $z'$ does not stabilise $B(\alpha,\beta)$.  

We claim that $\<x,(y_{1}^{g_k})^{z'}\>$ neither stabilises a block, nor fixes a point. As we noted above, $z' \cdot B(\alpha,\beta) \neq B(\alpha,\beta)$.  But \[y_1^{g_k}\cdot(z'\cdot\alpha)=z'x^ky_1x^{-k}z'^{-1}\cdot (z'\cdot\alpha)=z'\cdot\alpha;\] hence, $z'\cdot\alpha$ is one of the fixed points of $y_1^{g_k}$, which we assume stabilises $B(\alpha,\beta)$. We conclude that $z'\cdot\alpha=z'\cdot B(\alpha,\beta)\cap B(\alpha,\beta)$. By our choice of $k$, we have $z'\cdot\alpha\notin\{\alpha,\beta\}$.  Then the claim follows (because the only fixed block of $(y_1^{g_k})^{z'}$ is $z'\cdot B(\alpha,\beta)\neq B(\alpha,\beta)$).

If $\left<x,z'\right>$ is a structural subgroup, we can then apply the same arguments as for $\left<x,y\right>$. 

\noindent \textbf{Step 3:} $y$ does not fix either $\alpha$ or $\beta$, but stabilises $B(\alpha,\beta)$.\\
Since $z$ does not stabilise $B(\alpha,\beta)$ neither does $yz$. So we replace $y$ by $yz$, then either $\langle x,yz \rangle$ is not a structural subgroup or $yz$ fixes either $\alpha$ or $\beta$ and we can repeat steps 1 and 2.
\end{proof}

\begin{prop}\label{prop:eliminating structural subgroups - 2 fixed point case}
Let $(x,y,z) \in V_3(G)$ where $|x| \neq 2$ and where $x$ fixes exactly two points. Then we may connect $(x,y,z)$ to a generating $3$-tuple $(x,y',z')$ where $\left<x, x^{y'}\right>, \left<x,x^{z'}\right>$ are not structural subgroups. 
\end{prop}

\begin{proof}
We assume that $x$ fixes the points $\alpha$ and $\beta$ and stabilises the block $B = B(\alpha, \beta)$.  \Cref{prop:remove structural subgroups 1} implies that we may assume that $\left<x,y\right>$ and $\left<x,z\right>$ are both not structural subgroups where $y$ and $z$ are not involutions. In particular, $y,z \notin N_{G}(\left<x\right>)$. 

We will proceed in two stages first proving we can connect $(x,y,z)$ to a generating $3$-tuple $(x,y',z)$ where $\left<x, x^{y'}\right>$ is not a structural subgroup. If $\left<x,x^y\right>$ is not a structural subgroup, then there is nothing to prove. Otherwise, we note that $x^y \in G_{y \cdot \alpha, y \cdot \beta}$ and that $\left<x, x^y\right>$ is a subgroup of $G_\alpha$, $G_\beta$, or $G_B$. If $x^y\cdot B=B$, then $x$ stabilises $y^{-1}\cdot B\neq B$, a contradiction to Lemma \ref{lem:no-blocks}. Hence, we may assume that $\left<x,x^y\right> \leq G_\alpha$. Since $\left<x,y\right>$ is not a structural group, we have $y^{-1} \cdot \alpha \neq \alpha$.  Hence, we must have $y^{-1} \cdot \alpha = \beta$. Given that $\left<x,z\right>$ is not a structural subgroup, we must have $z \cdot \beta \neq \beta$. If $z \cdot \alpha = \beta$ and $z^{-1} \cdot  \alpha = \beta$, then $z$ stabilises $B$ and $\<x,z\>$ is structural, a contradiction. Thus, there exists a $z^\prime \in \{z,z^{-1}\}$ such that $z^\prime \cdot \beta \notin\{\alpha,\beta\}$ We also note that $z \cdot B \neq B$.

Let $g_i = z^\prime x^i y^{-1}$ for $0 \leq i \leq 3$. Then $g_i\cdot\alpha=z'\cdot\beta\notin\{\alpha,\beta\}$. We claim that $g_i \cdot \beta \neq g_j \cdot \beta$ for $0 \leq i < j \leq 3$. Suppose otherwise. We then have 
$x^{j-i} y^{-1} \cdot \beta = y^{-1} \beta$. It follows that $y^{-1} \cdot \beta \in \{\alpha, \beta\}$ since $x^{j-i}$ only fixes $\alpha$ and $\beta$. Given that $y^{-1} \cdot \beta \neq \beta$, we have $y^{-1} \cdot \beta = \alpha$. Since we assumed $y^{-1} \cdot \alpha = \beta$, this implies that $y\cdot B=B$, a contradiction. 

Therefore, $g_0 \cdot \beta, g_1 \cdot \beta, g_2 \cdot \beta, g_3\cdot\beta$ are all distinct. Hence, $g_k \cdot \beta \notin \{\alpha, \beta\}$ for at least two values of $k \in \{0,1,2,3\}$. Similarly, for at most one value of $k\in\{0,1,2,3\}$ we have $g_k\cdot B=B$. Hence, we see that there exists at least one $k\in\{0,1,2,3\}$ such that $g_k\cdot \beta\notin\{\alpha,\beta\}$ and $g_k\cdot B\neq B$. Since also $g_k\cdot \alpha\notin\{\alpha,\beta\}$, we conclude that $\<x,x^{g_k}\>$ is not structural.

We then connect
\[
(x,y,z) \to (x, y^{-1}, z) \to (x, x^ky^{-1},z) \to (x, z^\prime x^ky^{-1},z) = (x,g_k,z).
\]
Setting $y^\prime = g_k$, we may apply the same arguments to $(x,z,y')$ and $z$ to obtain a generating $3$-tuple $(x, y^\prime, z^\prime)$ where $\left<x,x^{z^\prime}\right>$ is not a structural subgroup.
\end{proof}

It is clear from the proofs that both \Cref{prop:remove structural subgroups 1} and \Cref{prop:eliminating structural subgroups - 2 fixed point case} also hold for generating $k$-tuples, with $k\geq 4$.

\subsection{$x$ fixes no points}
From now on, we again consider vectors in $V_5(G)$, that is, of length $5$. This allows us to again use \Cref{smallthenredundant} and \Cref{lem:eliminating small subgroups} in the course of our proofs.

\begin{prop}\label{prop:eliminating structural subgroups - fixed point free case}
Let $(x,y,z,v,w) \in V_5(G)$ be minimal and not connected to a redundant vector, where $x$ does not fix any point (in particular $|x| \neq 2$). Then either all subgroups $\<x,y,z,v\>$, $\<x,y,z,w\>$, $\<x,y,v,w\>$, $\<x,z,v,w\>$, $\<x,x^y\>$, $\<x,x^z\>$, $\<x,x^v\>$ and $\<x,x^w\>$ are non-structural, or  we can connect $(x,y,z,v,w)$ to a $5$-tuple containing a member of $G$ fixing exactly two points of $U_R(q)$.
\end{prop}

\begin{proof}
 By \Cref{lem:normalizer_involution}, we may assume that $y,z,v,w \notin N_G(\left<x\right>)$ and that $|y|,|z|, |v|, |w| \neq 2$. If $x$ does not fix any point, then, \hvm{by \Cref{lem:order-fix}}, its order divides either $\frac12(q+1)$, $q+\sqrt{3q}+1$, or $q-\sqrt{3q}+1$. In the latter two cases, $x$ is not contained in any structural subgroup, so that all of $\<x,x^y\>$, $\<x,x^z\>$, $\<x,x^v\>$, $\<x,x^w\>$, $\<x,y,z,v\>$, $\<x,y,z,w\>$, $\<x,y,v,w\>$ and $\<x,z,v,w\>$ are not structural by the minimality assumption. 

Hence, we may assume that $|x|$ divides $\frac12(q+1)$. 
We claim that we may assume that $\<y,z,v,w\>$ is a subfield subgroup. Indeed, suppose first that $\<y,z,v,w\>$ fixes some point $\alpha$. Then we can connect $(x,y,z,v,w)$ with $(x,xyx^{-1},z,u,v)$, and the subgroup $H:=\<xyx^{-1},z,v,w\>$ does not stabilise $\alpha$ or any other point of $\Omega$ (since all of $y,z,v,w$ fix exactly one point; if they fixed two points, then the assertion would trivially follow). 

Suppose now that $\<y,z,v,w\>$ stabilises a block $C$. Then we obtain that $\<y,z,v,w\>\leq \PSL_2(q)\times 2$. Since, by \cite{Wim:99}, all soluble subgroups of $\PSL_2(q)$ are $2$-generated,  we may, by \cite{DW}, assume that $\<y,z,v,w\>\in\{\PSL_2(q_0),\PSL_2(q_0)\times 2\}$, for some $3$-power $q_0$ with $q_0-1$ a divisor of $q-1$. If $\<y,z,v,w\>\cong\PSL_2(q_0)$, then, by \cite{G08}, we can connect $(y,z,v,w)$ to a redundant quadruple. We now show that $\<y,z,v,w\>\cong\PSL_2(q_0)\times 2$ leads to our assertion.  

 Let $\eta$ be  the unique involution fixing $C$ pointwise. We can write $(y,z,v,w)$ as \[(\eta^{\epsilon_y}y_0, \eta^{\epsilon_z}z_0,\eta^{\epsilon_v}v_0,\eta^{\epsilon_w}w_0),\] where $\epsilon_y,\ldots,\epsilon_w\in\{0,1\}$ and $y_0,\ldots,w_0\in\PSL_2(q_0)$. It follows from the main result of \cite{G08} that we can connect $(y_0,z_0,v_0,w_0)$ to $(y_0',z_0',v_0',w_0')$ with Nielsen moves, where $y_0'$ fixes exactly two points of $C$. The same Nielsen moves take $(y,z,v,w)$ to $(y^*,z^*,v^*,w^*)$, where $y_0'(y^*)^{-1},\ldots,w_0'(w^*)^{-1}\in\{1,\eta\}$. But then clearly, $y^*$ fixes exactly two points and the assertion follows. 

So we have shown that $\<y,z,v,w\>$ is not a structural subgroup. Since by \Cref{smallthenredundant} it is not a small subgroup either the claim follows.
 
 So, let $\<y,z,v,w\>\cong{^2\mathsf{G}}_2(q_1):=H$, for some appropriate power $q_1$ of $3$. \Cref{lem:eliminating small subgroups} immediately yields $q_1>3$. Let $b$ be the size of a block stabiliser in $H$. Then we claim that for at most $q_1b$ members $h$ of $H$, $y$ maps the block $h\cdot B$ to itself. Indeed, if $(yh)\cdot B=h\cdot B$ and $(yh')\cdot B=h'\cdot B$, then $yh'h^{-1}\cdot B=h'h^{-1}\cdot B$. Hence, $y$ stabilises $h'h^{-1}\cdot B$. The claim follows from the fact that $y$ stabilises at most $q_1$ blocks. Similarly, for at most $q_1b$ members $h$ of $H$, $y$ maps the block $h\cdot B$ to $B'$ and  for at most $q_1b$ members $h$ of $H$, $y$ maps the block $h\cdot B$ to $B''$. Similar statements hold for the blocks $B'$ and $B''$ and for $z,v$ and $w$. Since $36q_1b<|{^2\mathsf{G}}_2(q_1)|$, there exists a member $g\in H$ such that all of $y,z,v,w$ map $\{g\cdot B,g\cdot B',g\cdot B''\}$ to a disjoint set of blocks. It follows that $\<x^g,(x^g)^y\>$ is not structural, and the same thing holds for $\<x^g,(x^g)^z\>$, $\<x^g,(x^g)^v\>$ and $\<x^g,(x^g)^y\>$. The assertion now follows. 

The proof is complete.
\end{proof}

\section{Connecting to a redundant vector}\label{sec:redundant-connect}

In this section, we connect a generating $5$-tuple $(x,y,z,v,w)$ where $\left<x,x^y\right>$, $\left<x,x^z\right>$, $\left<x,x^v\right>$, and $\left<x,x^w\right>$ are non-structural subgroups to a redundant $5$-tuple. \Cref{unipotentconnectsto2or0} implies that we may assume that there exists an element of $(x,y,z,v,w)$ which fixes either no points or exactly two points. Thus, up to Nielsen moves, we may assume that $x$ is that element. By \Cref{prop:eliminating structural subgroups - 2 fixed point case} and \Cref{prop:eliminating structural subgroups - fixed point free case}, we may assume that $L_1 = \left<x,x^y\right>$, $L_2 =\left<x,x^z\right>$, $L_3 = \left<x,x^v\right>$, and $L_4 = \left<x,x^w\right>$ are not structural subgroups. Let $K_1 = \left<x,z,v,w\right>$, $K_2 =  \left<x,y,v,w\right>$, $K_3 = \left<x,y,z,w\right>$, and $K_4 = \left<x,y,z,v\right>$. We then see that $K_i$ are not structural subgroups since they contain non-structural subgroups for $1 \leq i \leq 4$.

\Cref{lem:eliminating small subgroups} implies that we may assume that $K_i$ are subfield subgroups over fields of size at least $27$ for $1 \leq i \leq 4$. In particular, we have $K_i \cong \: ^2\mathsf{G}_2(q_i)$ where $q_1 \geq q_2 \geq q_3 \geq q_4\geq 27$. 

\begin{prop}\label{prop:maximing |K1| + |K2|}
Using the above notation and assumptions, if $q_1 > q_4$, then one of the following holds:
\begin{enumerate}
\item $(x,y,z,v,w)$ is connected to a redundant $5$ tuple.
\item $(x,y,z,v,w)$ is connected to a generating $5$-tuple $(x,y',z,v,w)$ such that $\left<x,y'\right>$ is a subfield subgroup, where
\[
 |\left<x,y',z,v\right>| + |\left<x,y',z,w\right>| + |\left<x,y',v,w\right> | + |\left<x,z,v,w\right>| > \sum_{i=1}^4|K_i|.
\]
\end{enumerate}
\end{prop}
\begin{proof}
We claim that we can connect $(x,y,z,v,w)$ to $(x,y',z,v,w)$ such that 
\[\begin{cases}
 |C_{\left<x,z,v,w\right>}(x)| \leq |C_{\left<x,y',z,w\right>}(x)|, \\
 |C_{\left<x,z,v,w\right>}(x)| \leq  |C_{\left<x,y',z,v\right>}(x)|,
\mbox{ and }\\
|C_{\left<x,z,v,w\right>}(x)| \leq  |C_{\left<x,y',w,v\right>}(x)|.\end{cases}
\]
By \Cref{centralisers} $C_{\left<x,z,v,w\right>}(x)$ is one of the following: \begin{eqnarray*}
C_{\left<x,z,v,w\right>}(x) =
\begin{cases}
C_{{q_1}-1} & |x| \mid \frac{q-1}{2} \\
C_2^2 \times C_{\frac{q_1+1}{4}} & |x| \mid \frac{q+1}{4} \\
C_{q_1\pm\sqrt{3q_1}+1} & x \in M_i, i = \pm 1\\
P'\in \Syl_3({^2\mathsf{G}}_2(q_1)) & x \in Z(P), P \in \Syl_3(G)\\
[P',P'] \times \left<\eta \right>, \: |\eta| = 2 & x\in [P,P] \setminus Z(P), P \in \Syl_3(G)\\
Z(P') \cdot \left<x\right> & x \in P \setminus [P,P], P \in \Syl_3(G)\\
\left<x\right> \times \PSL_2(q_1) & |x|=2
\end{cases}
\end{eqnarray*}
Therefore, we have a few cases. The first case takes the first and the third possibility above together and corresponds to $C_{\left<x,z,v,w\right>}(x)$ being cyclic. The second case is when $C_{\left<x,z,v,w\right>}(x) \cong C_2^2 \times C_{\frac{q_1+1}{4}}$. The last four possibilities above are when $x$ is either an involution or has order $3,6$ or $9$. However, \Cref{unipotentconnectsto2or0} implies that we may assume that $x$ does not have order $2,3,6,$ or $9$. 

\noindent \textbf{Case 1:} $C_{\left<x,z,v,w\right>}(x)$ is cyclic.\\
In this case, we have $C_{\left<x,z,v,w\right>}(x) \in \{C_{q_1-1}, M_{-1}, M_1\}$. Since $x$ and $x^y$ are conjugate and $\left<x,x^y\right>$ is not a structural group, by \Cref{lem:conjugate in smaller ree group} and our assumption on $K_i$ being subfield subgroups there exists an element $d \in \left<x,x^y\right>$ such that $x^y = x^d$. Let $u =yd^{-1}$, which gives $u \in C_{\left<x,y\right>}(x)$. Writing  $C_{\left<x,z,v,w\right>}(x) = \left<c\right>$, we claim that $\left<u,c\right>$ is cyclic, and we proceed based on the isomorphism type of $C_{\left<x,z,v,w\right>}.$ Since $x$ does not have order $2,3,6$ or $9$, we have $C_{G}(x) \in \{C_{q-1}, M_{-1}, M_1, C_2^2 \times C_{\frac{q+1}{4}}\}$. 

If $C_{\left<x,z,v,w\right>}(x) = C_{{q_1}-1}$, we then have $|x| \mid \frac{q_1-1}{2}$, which implies $2|x| \mid q_1-1$. Since $q_1 = 3^{2m+1}$ and $q=3^{2n+1}$, where $2m+1 \mid 2n+1$, we have $q_1 -1  \mid q-1$. Therefore, $2|x| \mid q-1$; hence, $|x| \mid \frac{q-1}{2}$. Hence, $C_{G}(x) \cong C_{q-1}$, which is cyclic. Since $\left<u,c\right> \leq C_{G}(x)$, it follows that $\left<c,u\right>$ is cyclic. 

Now suppose that $C_{\left<x,z,v,w\right>}(x)$ is cyclic of order $q_1 \pm \sqrt{3q_1} + 1$.  In this case the order of $x$ divides $q\pm\sqrt{3q}+1$ and not $q+1$, so $C_G(x)$ is again cyclic, and so is $\<c,u\>$.  

  \Cref{lem:cyclic subgroup reduction} then implies there exists a natural number $m$ such that $\left<c,u\right> =  \left<uc^m\right>$. Since $c^m \in \left<x,z,v,w\right>$, we may connect $(x,y,z,v,w)$ to $(x,yc^m,z,v,w)$. Since
\[
d \in \left<x, x^{y}\right> = \left<x,x^{yc^m} \right> \leq \left<x,yc^m, z,w\right> \cap \left<x,yc^m,z,a\right> \cap \left<x,yc^m,w,v\right>,
\]
we see that $\left<x,yc^m, z, w\right>$, $\left<x,yc^m,z,v\right>$ and $\left<x,yc^m,w,v\right>$ are subgroups that contain $d$ which are not structural. Thus,
\[
d^{-1}yc^m = uy^{-1}yc^m  = uc^m \in \left<x,yc^m, z,w\right> \cap \left<x,yc^m,z,v\right> \cap \left<x,yc^m,w,v\right>.
\]
Given that $C_{\left<x,z,v,w\right>}(x) = \left<c\right> \in \left<uc^m\right>$, we have
\[
C_{\left<x,z,v,w\right>}(x) \leq C_{\left<x,yc^m, z,w\right>}(x) \cap C_{\left<x,yc^m,z,v\right>}(x) \cap C_{\left<x,yc^m,w,v\right>}(x).
\]
Thus, we have $|C_{\left<x,z,v,w\right>}(x)| \leq |C_{\left<x,yc^m,z,v\right>}(x)|$. Similarly, 
\[
|C_{\left<x,z,v,w\right>}(x)| \leq |C_{\left<x,yc^m,w,v\right>}(x)| \quad \text{and} \quad |C_{\left<x,z,v,w\right>}(x)| \leq |C_{\left<x,yc^m,z,w\right>}(x)|.
\]
\Cref{lem:centralizer ree subgroup inequality} implies 
\[
|\left<x,z,v,w\right>| \leq |\left<x,yc^m,z,v\right>|, \quad |\left<x,z,v,w\right>| \leq |\left<x,yc^m,w,v\right>|,
\]
and
\[
|\left<x,z,v,w\right>| \leq |\left<x,yc^m,z,w\right>|.
\]

\noindent \textbf{Case 2:} $C_{\left<x,z,v,w\right>}(x)  = C_2^2 \times C_{\frac{q_1+1}{4}}$\\ 
This case happens when $|x| \mid \frac{q_1 + 1}{4}$. Since $\left<x,x^y\right>$ is not a structural group, by \Cref{lem:conjugate in smaller ree group} and our assumption that the $K_i$ are subfield subgroups there exists an element $d \in \left<x,x^y\right>$ such that $x^y = x^d$. Let $u =yd^{-1}$, which gives $u \in C_{\left<x,y\right>}(x)$. Let $c\in\<x,z,v,w\>$ be an order $\frac14(q_1+1)$ element such that $C_{\<x,z,v,w\>}(x)= C_2^2\times\<c\>$, where $C_2^2$ is generated by some two involution $\eta_1,\eta_2$. Since $|x|$ divides $\frac12(q_1+1)$, it also divides $\frac12(q+1)$, and hence, there exists an order $\frac14(q+1)$ element $b\in G$ such that $C_G(x)= \<\eta_1,\eta_2\>\times\<b\>$. Then we have $c=b^k$, for some $k\in\mathbb{N}$, and $u=\eta^\epsilon b^j$, for some $j\in\mathbb{N}$, $\eta\in \<\eta_1,\eta_2\>$ and $\epsilon\in\{0,1\}$. It follows that $\<c,u\>\leq\<\eta\>\times\<b\>$, which is cyclic. 

  \Cref{lem:cyclic subgroup reduction} then implies that there exists a natural number $m$ such that $\left<c,u\right> =  \left<uc^m\right>$, and the rest of the proof of Case 1 can be copied. 
\end{proof}

In case $q_1 = q_4$ we obtain the following proposition.
\begin{prop}\label{prop: q_1 = q_4 finish}
With respect to the above notation, if $q_1 =q_4$, then $(x,y,z,v,w)$ is connected to a redundant $5$-tuple
\end{prop}
\begin{proof}
The groups $K_1, K_2, K_3,$ and $K_4$ are all isomorphic to $^2\mathsf{G}_2(q_1)$. Let $H=\left<x,v,w\right> \leq K_1 \cap K_2$.  Since $L_3\leq H$ and $L_3$ is not structural, $H$ is either small or a subfield subgroup. Suppose $H$ is small. Then it is solvable and generated by two elements. By \cite{DW}, $(x,v,w)$ is connected to a redundant triple, and we are done in this case. Hence, we may assume that $H$ is a subfield subgroup, and hence self-normalising.

Since $K_1$ and $K_2$ are subfield subgroups of the same order, there exists an element $g \in G$ so that $K_1^g = K_2$. Thus, since $H\leq K_1$, we have $H^g\leq K_2$. However, $H$ is also a subgroup of $K_2$. Therefore, $H$ and $H^g$ are conjugate in $K_2$. Hence, as $K_2$ is a subfield subgroup, there exists an element $h \in K_2$ such that $H^g = H^h$. So, $h^{-1}g \in N_{G}(H)=H \leq K_2$. Consequently, $g \in K_2$. We now have $K_1 = K_2^{g^{-1}} =K_2$. Since $K_1 = \left<x,z,v,w\right>$ and $K_2 = \left<x,y,v,w\right>$, we have $G = \left<K_1, K_2\right>$. The proposition follows. 
\end{proof}

We now come to the proof of our Main Theorem.

 \begin{proof}
Let $(x,y,z,v,w) \in V_5(G)$ be a generating $5$-tuple. Using \Cref{unipotentconnectsto2or0} implies that we may assume that $x$ either fixes exactly two points or fixes no points. Moreover, using \Cref{prop:eliminating structural subgroups - 2 fixed point case} and \Cref{prop:eliminating structural subgroups - fixed point free case} we may also assume, with previous notation, that $L_i$ and $K_i$ are not structural subgroups, $i=1,2,3,4$.
We note that if any of the subgroups $K_1, K_2, K_3, K_4$ are small, we then have that $(x,y,z,v,w)$ is connected to a redundant $5$-tuple, by \Cref{lem:eliminating small subgroups}. Therefore, we may assume that $K_1, K_2, K_3$, and $K_4$ are subfield subgroups. 
We may apply \Cref{prop:maximing |K1| + |K2|}, possibly interchanging the roles of $y,z,v$ and $w$, finitely many times until we either can connect $(x,y,z,v,w)$ to a redundant $5$-tuple, or $|K_1| = |K_2| = |K_3| = |K_4|$. \Cref{prop: q_1 = q_4 finish} implies that $(x,y,z,v,w)$ is connected to a redundant $5$-tuple, as desired.
\end{proof}

\appendix

\section{Matrices for the point stabilisers}\label{appendixA}
In this appendix we review the explicit construction of $^2\mathsf{G}_2(q)$ given in \cite[Section~9.2.4]{TTM07}. It is used in the present paper to prove various claims about the action of $^2\mathsf{G}_2(q)$ on $U_R(q)$, in particular actions on blocks of this unital.

We start with defining an explicit model of $U_R(q)$. Recall $q=3^{2e+1}$ and $\theta:\mathbb{F}_q\to\mathbb{F}_q:x\mapsto x^{3^{e+1}}$. We then define $U_R(q)$ to be the set of points of the projective $6$-space $\PG(6,q)$ with coordinates given by $(1,0,0,0,0,0,0)$ and \[(f_1(a,a',a''),-a',-a,-a'',1,f_2(a,a',a''),f_3(a,a',a''))=:p(a,a',a''),\] $a,a',a''\in\mathbb{F}_q$ and 
\[\begin{cases}
f_1(a,a',a'')=-a^{2\theta+4}-aa''^\theta+a^{\theta+1}a'^\theta+a''^2+a'^{\theta+1}-a'a^{\theta+3}-a^2a'^2,\\
f_2(a,a',a'')=-a^{\theta+3}+a'^\theta-aa''+a^2a',\\
f_3(a,a',a'')=-a^{2\theta+3}-a''^\theta+a^\theta a'^\theta+a'a''+aa'^2.
\end{cases}\]

We set $\infty=(1,0,0,0,0,0,0)$ and $O=(0,0,0,0,1,0,0)$. Then the Sylow $3$-subgroup $P_\infty$ of $^2\mathsf{G}_2(q)$ fixing $\infty$ is given by the set of linear transformations $(a,a',a'')_\infty$ (acting on the left) with corresponding matrices
\[\begin{pmatrix}
1 & p & q & a'' & f_1(a,a',a'') & a'-a^{\theta+1} & a \\
0 & 1 & a^\theta & 0 & -a' & 0 & 0\\
0 & 0 & 1 & 0 & -a & 0 & 0 \\
0 & -a & a'-a^{\theta+1} & 1 & -a''& 0 & 0 \\
0 & 0 & 0 & 0 & 1 & 0 & 0\\
0 & a^2 & r & a & f_2(a,a',a'') & 1 & 0 \\
0 & -a''-aa' & s & -a' & f_3(a,a',a'') & -a^\theta & 1
\end{pmatrix},\]
where 
\[\begin{cases}
p=a^{\theta+3}-a'^\theta-aa''-a^2a',\\
q=-a^{2\theta+3}+a''^\theta + a^\theta a'^\theta + a'a''-aa'^2-a^{\theta+2}a'-a^{\theta+1}a'',\\
r=a^{\theta+2}+a''-aa',\\
s=-a^\theta a''+a'^2-a^{\theta+1}a'.
\end{cases}
\]
The Sylow $3$-subgroup $P_O$ of $^2\mathsf{G}_2(q)$ fixing $O$ is given by the set of linear transformations $(a,a',a'')_O$ (again acting on the left) with corresponding matrices
\[\begin{pmatrix}
 1 & 0 & 0 & 0 & 0 & 0 & 0 \\
f_2(a,a',a'') & 1 & 0 & -a & 0 & a^2 & r  \\
  f_3(a,a',a'') & -a^\theta & 1 & a' & 0 & -a''-aa' & s \\
 a''& 0 & 0 & 1 & 0 & a & a^{\theta+1}-a'  \\
 f_1(a,a',a'') & a'-a^{\theta+1} & a & -a'' &1 & p & q  \\
 -a' & 0 & 0 & 0 & 0 & 1 & a^\theta \\
  -a & 0 & 0 & 0 & 0 & 0 & 1 
\end{pmatrix}.\]
Its action on $U_R(q)$ can best be seen using an alternative description of $U_R(q)$ (see also \cite[Section~9.2.4]{TTM07}), namely $U_R(q)$ consists of the point $(0,0,0,0,1,0,0)$ and the points  \[q(a,a',a''):=(1,f_2(a,a',a''),f_3(a,a',a''),a'',f_1(a,a',a''),-a',-a),\] $a,a',a''\in\mathbb{F}_q$. These points are the images of $(1,0,0,0,0,0)$ under the maps $(a,a',a'')_O$, as is readily checked. 

Now we note that the matrices displayed above all have determinant equal to $1$. Since $^2\mathsf{G}_2(q)\cong\<P_\infty,P_O\>$, we obtain $^2\mathsf{G}_2(q)\leq\SL_7(q)$. Also, since $7$ does not divide $q-1$, every element of $^2\mathsf{G}_2(q)$ is represented by exactly one matrix, and so we can define the trace of an element as the trace of its matrix. In view of the matrices given above, this immediately leads to the following observation.

\begin{obs}\label{tracelemma}
If $x\in{^2\mathsf{G}}_2(q)$ has order $3$ or $9$, then (the corresponding matrix of) $x$ (in $\SL_7(q)$) has trace~$1$. 
If $x\in{^2\mathsf{G}}_2(q)$ has order $2$ or $6$, then (the corresponding matrix of) $x$ (in $\SL_7(q)$) has trace~$-1$. 
\end{obs}

One verifies easily the following multiplication law: $(a,a',a'')_\infty(b,b',b'')_\infty=$
\[(a+b,a'+b'+ a^\theta b,a''+b''-ab'+a'b-a^{\theta+1}b)_\infty.\] Note that this differs slightly from the expression in  \cite[Section~9.2.4]{TTM07} as we consider here action on the left, whereas the action is on the right in  \cite[Section~9.2.4]{TTM07}. We deduce \[(a,a',a'')_\infty\cdot p(b,b',b'')=p(a+b,a'+b'+ a^\theta b,a''+b''-ab'+a'b-a^{\theta+1}b).\] The inverse is given by \[(a,a',a'')_\infty^{-1}=(-a,-a'+a^{\theta+1},-a'')_\infty.\] It is easy to check now that \[\begin{cases}Z(P_\infty)=\{(0,0,a'')_\infty\mid a''\in\mathbb{F}_q\},\\ [P_\infty,P_\infty]=\{(0,a',a'')_\infty\mid a',a''\in\mathbb{F}_q\}\\ \{x^3\mid x\in P_\infty\} = Z(P_\infty).\end{cases}\]

The two-point stabiliser $G_{\infty,O}$ is given  the following action on $U_R(q)$. 
\[h(t)\cdot p(b,b',b'')=p(bt^\theta,b't^{\theta+1},b''t^{\theta+2}), t,b,b',b''\in\mathbb{F}_q.\] We find that $h(-1)$ is an involution fixing the point set $\{\infty\}\cup\{p(0,b',0)\mid b'\in\mathbb{F}_q\}$, which is --- by definition --- the block $B(\infty,O)$. The other blocks through $\infty$ can be obtained by the image of $B(\infty,O)$ under the elements $(a,0,a'')_\infty$, $a,a''\in\mathbb{F}_q$ (since $(0,a',0)_\infty$ stabilises $B(\infty,O)$). A generic block like that looks like \[\{p(a,b',a''-ab')\mid b'\in\mathbb{F}_q\}=:B_{a,a''}.\]
The following is immediate.
\begin{obs}\label{obs1}
Every nontrivial member of $G_{\infty,O}$ fixes exactly two points ($\infty$ and $O$), except $h(-1)$, which fixes $B(\infty,O)$ pointwise, but no point off that block. Moreover, $G_{\infty,O}$ acts freely on $U_R(q)\setminus B(\infty,O)$.  
\end{obs}
The second assertion of the previous observation follows from the fact that, if $p(bt^\theta,b't^{\theta+1},b''t^{\theta+2})=p(bu^\theta,b'u^{\theta+1},b''u^{\theta+2})$, for arbitrary $b,b',b'',t,u\in\mathbb{F}_q$, then either $b=b''=0$, or $t=u$ (noting that $t\mapsto t^{\theta+2}$ is injective).

Using the explicit form of the multiplication, inverse and action of the two-point stabiliser, one calculates that the elements $(1,0,0)_\infty$, $(1,1,0)_\infty$ and $(1,-1,0)_\infty$ are representatives of the three different conjugacy classes of elements of order $9$ of $P_\infty$, whereas $(0,0,1)_\infty$, $(0,1,0)_\infty$ and $(0,-1,0)_\infty$ are representatives of the three different conjugacy classes of elements of order $3$ of $P_\infty$. Also, \[p(b,b',b'')\mapsto p(-b,b'+1,-b''-b) \mbox{ and }p(b,b',b'')\mapsto p(-b,b'-1,-b''+b)\] are representatives of the two conjugacy classes of elements of order $6$. 

Finally we make the following observation.

\begin{obs}\label{calculation}
The point $p(0,a',a'')$, with $a''\neq 0$, coincides with the point $q(0,b',b'')$, $b''\neq 0$, if, and only if, $b'=b''^{\theta-1}=a''^{1-\theta}=a'^{-1}$.
\end{obs}

\begin{proof}We write the points $p(0,a',a'')$ and $q(0,b',b'')$ in projective coordinates and obtain
\[\begin{cases} p(0,a',a'')=(a''^2+a'^{\theta+1},-a',0,-a'',1,a'^\theta,-a''^\theta+a'a'')=:p,\\ q(0,b',b'')= (1,b'^\theta,-b''^\theta+b'b'',b'',b''^2+b'^{\theta+1},-b',0)=:q.\end{cases}\]
If $p=q$, then clearly $a'a''=a''^\theta$ and $b'b''=b''^\theta$. Hence, $a'=a''^{\theta-1}$ and $b'=b''^{\theta-1}$. Since projective coordinates are determined up to a non-zero scalar factor, we can rewrite the projective coordinates of $p$ and $q$ as follows:
\[\begin{cases}
p=(-a''^2,-a''^{\theta-1},0,-a'',1,a''^{3-\theta},0),\\
q=(-b''^{-2},-b''^{1-\theta},0,-b''^{-1},1,b''^{\theta-3},0),
\end{cases}\]
from which the assertion readily follows, noting that the converse is easy. 
\end{proof}\color{black}

\end{document}